\input amstex
\documentstyle{amsppt}
%
\catcode`@=11
\redefine\output@{%
  \def\break{\penalty-\@M}\let\par\endgraf
  \ifodd\pageno\global\hoffset=105pt\else\global\hoffset=8pt\fi  
  \shipout\vbox{%
    \ifplain@
      \let\makeheadline\relax \let\makefootline\relax
    \else
      \iffirstpage@ \global\firstpage@false
        \let\rightheadline\frheadline
        \let\leftheadline\flheadline
      \else
        \ifrunheads@ 
        \else \let\makeheadline\relax
        \fi
      \fi
    \fi
    \makeheadline \pagebody \makefootline}%
  \advancepageno \ifnum\outputpenalty>-\@MM\else\dosupereject\fi
}
\def\Beta{\mathchar"0\hexnumber@\rmfam 42}
\catcode`\@=\active
\nopagenumbers
\chardef\textvolna='176
\def\negskp{\hskip -2pt}
\def\Ker{\operatorname{Ker}}

\def\Sym{\operatorname{Sym}}
\chardef\degree="5E
\def\blue#1{#1}

\catcode`#=11\def\diez{#}\catcode`#=6
\catcode`&=11\catcode`&=4
\catcode`_=11\def\podcherkivanie{_}\catcode`_=8
\catcode`~=11\def\volna{~}\catcode`~=\active
\def\mycite#1{\cite{\blue{#1}}\immediate\special{ps:
     ShrHPSdict begin /ShrBORDERthickness 0 def}}
\def\myciterange#1#2#3#4{\cite{\blue{#2#3#4}}\immediate\special{ps:
     ShrHPSdict begin /ShrBORDERthickness 0 def}}
\def\mytag#1{%
    \tag#1}
\def\mythetag#1{\thetag{\blue{#1}}\immediate\special{ps:
     ShrHPSdict begin /ShrBORDERthickness 0 def}}
\def\myrefno#1{\no#1}
\def\myhref#1#2{\blue{#2}\immediate\special{ps:
     ShrHPSdict begin /ShrBORDERthickness 0 def}}
\def\myEarXivlink{\myhref{http://arXiv.org}{http:/\negskp/arXiv.org}}

\def\mytheorem#1{\csname proclaim\endcsname{Theorem #1}}
\def\mytheoremwithtitle#1#2{\csname proclaim\endcsname{Theorem #1#2}}
\def\mythetheorem#1{\blue{#1}\immediate\special{ps:
     ShrHPSdict begin /ShrBORDERthickness 0 def}}
\def\mylemma#1{\csname proclaim\endcsname{Lemma #1}}
\def\mylemmawithtitle#1#2{\csname proclaim\endcsname{Lemma #1#2}}

\def\mycorollary#1{\csname proclaim\endcsname{Corollary #1}}

\def\mydefinition#1{\definition{Definition #1}}
\def\mythedefinition#1{\blue{#1}\immediate\special{ps:
     ShrHPSdict begin /ShrBORDERthickness 0 def}}
\def\myconjecture#1{\csname proclaim\endcsname{Conjecture #1}}
\def\myconjecturewithtitle#1#2{\csname proclaim\endcsname{Conjecture #1#2}}

\def\myanchortext#1#2{#2}
\def\mytheanchortext#1#2{\blue{#2}\immediate\special{ps:
     ShrHPSdict begin /ShrBORDERthickness 0 def}}
\font\eightcyr=wncyr8
\pagewidth{360pt}
\pageheight{606pt}
\topmatter
\title
A biquadratic Diophantine equation associated with perfect cuboids.
\endtitle
\rightheadtext{A biquadratic Diophantine equation \dots}
\author
Ruslan Sharipov
\endauthor
\address Bashkir State University, 32 Zaki Validi street, 450074 Ufa, Russia
\endaddress
\email\myhref{mailto:r-sharipov\@mail.ru}{r-sharipov\@mail.ru}
\endemail
\abstract
    A perfect Euler cuboid is a rectangular parallelepiped with integer edges and 
integer face diagonals whose space diagonal is also integer. Such cuboids are not
yet discovered and their non-existence is also not proved. Perfect Euler cuboids
are described by a system of four Diophantine equation possessing a natural $S_3$
symmetry. Recently these equations were factorized with respect to this $S_3$ 
symmetry and the factor equations were derived. In the present paper the factor 
equations are transformed to $E$-form and then reduced to a single biquadratic 
equation. 
\endabstract
\subjclassyear{2000}
\subjclass 11D41, 11D72, 13A50, 13F20\endsubjclass
\endtopmatter
\TagsOnRight
\document

\head
1. Introduction.
\endhead
     Perfect cuboids are described by four polynomial equations 
$$
\xalignat 4
&\hskip -2em
p_{\kern 1pt 0}=0,&&p_{\kern 1pt 1}=0,&&p_{\kern 1pt 2}=0,&&p_{\kern 1pt 3}=0,
\quad
\mytag{1.1}
\endxalignat
$$
where $p_0$, $p_1$, $p_2$, $p_3$ are the following polynomials of seven variables:
$$
\xalignat 2
&\hskip -2em
p_{\kern 1pt 0}=x_1^2+x_2^2+x_3^2-L^2,
&&p_{\kern 1pt 1}=x_2^2+x_3^2-d_1^{\kern 1pt 2},\\
\vspace{-1.7ex}
\mytag{1.2}\\
\vspace{-1.7ex}
&\hskip -2em
p_{\kern 1pt 2}=x_3^2+x_1^2-d_2^{\kern 1pt 2},
&&p_{\kern 1pt 3}=x_1^2+x_2^2-d_3^{\kern 1pt 2}.
\endxalignat
$$
Here $x_1$, $x_2$, $x_3$ are edges of a cuboid, $d_1$, $d_2$, $d_3$ are its face 
diagonals, and $L$ is its space diagonal. Though the equations \mythetag{1.1}
with the polynomials \mythetag{1.2} look very simple, the search for perfect 
cuboids has the long history since 1719 (see \myciterange{1}{1}{--}{44}).\par
     Recently in \mycite{45} the symmetry approach to the equations \mythetag{1.1}
was initiated. It is based on the intrinsic $S_3$ symmetry of these equations. 
Let the permutation group $S_3$ act upon the variables $x_1$, $x_2$, $x_3$, $d_1$, 
$d_2$, $d_3$, $L$ according to the rules
$$
\xalignat 3
&\hskip -2em
\sigma(x_i)=x_{\sigma i},
&&\sigma(d_i)=d_{\sigma i},
&&\sigma(L)=L.
\mytag{1.3}
\endxalignat
$$ 
If the variables $x_1$, $x_2$, $x_3$ and $d_1$, $d_2$, $d_3$ are arranged into
the matrix
$$
\hskip -2em
M=\Vmatrix x_1 & x_2 &x_3\\
\vspace{1ex}
d_1 & d_2 & d_3\endVmatrix,
\mytag{1.4}
$$
then, according to \mythetag{1.3}, the group $S_3$ acts upon the matrix \mythetag{1.4} 
by permuting its columns. Applying the rules \mythetag{1.3} to the polynomials 
\mythetag{1.2}, we derive 
$$
\xalignat 2
&\hskip -2em
\sigma(p_i)=p_{\kern 0.5pt \sigma i},
&&\sigma(p_{\kern 1pt 0})=p_{\kern 1pt 0}.
\mytag{1.5}
\endxalignat
$$
The polynomials $p_{\kern 1pt 0}$, $p_{\kern 1pt 1}$, $p_{\kern 1pt 2}$, 
$p_{\kern 1pt 3}$ in \mythetag{1.2} and in \mythetag{1.5} are elements of the polynomial
ring $\Bbb Q[x_1,x_2,x_3,d_1,d_2,d_3,L]$. For the sake of brevity we denote this ring
through
$$
\hskip -2em
\Bbb Q[M,L]=\Bbb Q[x_1,x_2,x_3,d_1,d_2,d_3,L],
\mytag{1.6}
$$
where $M$ is the matrix given by the formula \mythetag{1.4}. 
\mydefinition{1.1} A polynomial $p\in\Bbb Q[M,L]$ is called multisymmetric if it
is invariant with respect to the action \mythetag{1.3} of the group $S_3$. 
\enddefinition
     General multisymmetric polynomials, which are also known as vector 
symmetric polynomials, diagonally symmetric polynomials, McMahon polynomials 
etc, were initially studied in \myciterange{46}{46}{--}{52} (see also later 
publications \myciterange{53}{53}{--}{66}).\par
     Multisymmetric polynomials constitute a subring in the ring \mythetag{1.6}.
We denote this subring through $\Sym\!\Bbb Q[M,L]$. The formulas \mythetag{1.5}
show that the polynomial $p_{\kern 1pt 0}$ belongs to the subring $\Sym\!\Bbb 
Q[M,L]$, i\.\,e\. it is multisymmetric, while the polynomials $p_{\kern 1pt 1}$, 
$p_{\kern 1pt 2}$, $p_{\kern 1pt 3}$ are not multisymmetric. Nevertheless, the 
system of equations \mythetag{1.1} in whole is invariant with respect to the 
action of the group $S_3$.\par
     The polynomials $p_{\kern 1pt 0}$, $p_{\kern 1pt 1}$, $p_{\kern 1pt 2}$, 
$p_{\kern 1pt 3}$ generate an ideal in the ring $\Bbb Q[M,L]$. In \mycite{67}
this ideal was called the {\it perfect cuboid ideal\/} and was denoted through
$$
\hskip -2em
I_{\text{PC}}=\bigl<p_{\kern 1pt 0},p_{\kern 1pt 1},p_{\kern 1pt 2},
p_{\kern 1pt 3}\bigr>.
\mytag{1.7}
$$
Each polynomial equation $p=0$ with $p\in I_{\text{PC}}$ follows from the 
equations \mythetag{1.1}. Therefore such an equation is called a {\it perfect
cuboid equation}.\par
     The symmetry approach to the equations \mythetag{1.1} initiated in 
\mycite{45} leads to studying the following ideal in the ring of multisymmetric 
polynomials $\Sym\!\Bbb Q[M,L]$: 
$$
\hskip -2em
I_{\text{PC\kern -0.7pt\_\kern 0.5pt sym}}=I_{\text{PC}}\cap\Sym\!\Bbb Q[M,L].
\quad
\mytag{1.8}
$$
\mydefinition{1.2} A polynomial equation of the form $p=0$ with $p\in 
I_{\text{PC\kern -0.7pt\_\kern 0.5pt sym}}$ is called an {\it $S_3$ factor 
equation} for the perfect cuboid equations \mythetag{1.1}. 
\enddefinition
     The ideal \mythetag{1.7} was initially studied in \mycite{68}. There it was
denoted through $I_{\text{sym}}$. However, in this paper we use the notation 
\mythetag{1.8} taken from \mycite{67}. In \mycite{68}, when studying the ideal 
$I_{\text{PC\kern -0.7pt\_\kern 0.5pt sym}}$, the following eight polynomials 
were introduced:
$$
\allowdisplaybreaks
\gather
\hskip -2em
\tilde p_{\kern 1pt 1}=p_{\kern 1pt 0}=x_1^2+x_2^2+x_3^2-L^2.
\mytag{1.9}\\
\vspace{2ex}
\hskip -2em
\aligned
\tilde p_{\kern 1pt 2}&=p_{\kern 1pt 1}+p_{\kern 1pt 2}+p_{\kern 1pt 3}
=(x_2^2+x_3^2-d_1^{\kern 1pt 2})\,+\\
&+\,(x_3^2+x_1^2-d_2^{\kern 1pt 2})+(x_1^2+x_2^2-d_3^{\kern 1pt 2}),
\endaligned\qquad\qquad
\mytag{1.10}\\
\vspace{2ex}
\aligned
\tilde p_{\kern 1pt 3}&=d_1\,p_{\kern 1pt 1}+d_2\,p_{\kern 1pt 2}+d_3
\,p_{\kern 1pt 3}=d_1\,(x_2^2+x_3^2-d_1^{\kern 1pt 2})\,+\\
&+\,d_2\,(x_3^2+x_1^2-d_2^{\kern 1pt 2})+d_3\,(x_1^2+x_2^2-d_3^{\kern 1pt 2}),
\endaligned\qquad\qquad
\mytag{1.11}\\
\vspace{2ex}
\aligned
\tilde p_{\kern 1pt 4}&=x_1\,p_{\kern 1pt 1}+x_2\,p_{\kern 1pt 2}+x_3
\,p_{\kern 1pt 3}=x_1\,(x_2^2+x_3^2-d_1^{\kern 1pt 2})\,+\\
&+\,x_2\,(x_3^2+x_1^2-d_2^{\kern 1pt 2})+x_3\,(x_1^2+x_2^2-d_3^{\kern 1pt 2}),
\endaligned\qquad\qquad
\mytag{1.12}\\
\vspace{2ex}
\aligned
\tilde p_{\kern 1pt 5}&=x_1\,d_1\,p_{\kern 1pt 1}+x_2\,d_2\,p_{\kern 1pt 2}
+x_3\,d_3\,p_{\kern 1pt 3}=x_1\,d_1\,(x_2^2+x_3^2-d_1^{\kern 1pt 2})\,+\\
&+\,x_2\,d_2\,(x_3^2+x_1^2-d_2^{\kern 1pt 2})+x_3\,d_3\,(x_1^2+x_2^2
-d_3^{\kern 1pt 2}),
\endaligned\qquad\qquad
\mytag{1.13}\\
\vspace{2ex}
\aligned
\tilde p_{\kern 1pt 6}&=x_1^2\,p_{\kern 1pt 1}+x_2^2\,p_{\kern 1pt 2}+x_3^2
\,p_{\kern 1pt 3}=x_1^2\,(x_2^2+x_3^2-d_1^{\kern 1pt 2})\,+\\
&+\,x_2^2\,(x_3^2+x_1^2-d_2^{\kern 1pt 2})+x_3^2\,(x_1^2+x_2^2
-d_3^{\kern 1pt 2}),
\endaligned\qquad\qquad
\mytag{1.14}\\
\vspace{2ex}
\aligned
\tilde p_{\kern 1pt 7}&=d_1^{\kern 1pt 2}\,p_{\kern 1pt 1}+d_2^{\kern 1pt 2}
\,p_{\kern 1pt 2}+d_3^{\kern 1pt 2}\,p_{\kern 1pt 3}=d_1^{\kern 1pt 2}
\,(x_2^2+x_3^2-d_1^{\kern 1pt 2})\,+\\
&+\,d_2^{\kern 1pt 2}\,(x_3^2+x_1^2-d_2^{\kern 1pt 2})+d_3^{\kern 1pt 2}
\,(x_1^2+x_2^2-d_3^{\kern 1pt 2}),
\endaligned\qquad\qquad
\mytag{1.15}\\
\vspace{2ex}
\aligned
\tilde p_{\kern 1pt 8}&=x_1^2\,d_1^{\kern 1pt 2}\,p_{\kern 1pt 1}+x_2^2
\,d_2^{\kern 1pt 2}\,p_{\kern 1pt 2}+x_3^2\,d_3^{\kern 1pt 2}
\,p_{\kern 1pt 3}=x_1^2\,d_1^{\kern 1pt 2}\,(x_2^2+x_3^2-d_1^{\kern 1pt 2})\,+\\
&+\,x_2^2\,d_2^{\kern 1pt 2}\,(x_3^2+x_1^2-d_2^{\kern 1pt 2})+x_3^2
\,d_3^{\kern 1pt 2}\,(x_1^2+x_2^2-d_3^{\kern 1pt 2}).
\endaligned\qquad\qquad
\mytag{1.16}
\endgather
$$
The polynomials $\tilde p_{\kern 1pt 0}$, $\tilde p_{\kern 1pt 1}$, 
$\tilde p_{\kern 1pt 2}$, $\tilde p_{\kern 1pt 3}$, $\tilde p_{\kern 1pt 4}$, 
$\tilde p_{\kern 1pt 5}$, $\tilde p_{\kern 1pt 6}$, $\tilde p_{\kern 1pt 7}$, 
$\tilde p_{\kern 1pt 8}$ in \mythetag{1.9}, \mythetag{1.10}, \mythetag{1.11}, 
\mythetag{1.12}, \mythetag{1.13}, \mythetag{1.14}, \mythetag{1.15}, 
\mythetag{1.16} are multisymmetric and they are generated by the polynomials 
$p_{\kern 1pt 0}$, $p_{\kern 1pt 1}$, $p_{\kern 1pt 2}$, $p_{\kern 1pt 3}$
from \mythetag{1.2}. For this reason they belong to the ideal \mythetag{1.8}. 
Moreover, there is the following theorem.
\mytheorem{1.1}The ideal $I_{\text{PC\kern -0.7pt\_\kern 0.5pt sym}}$ in
\mythetag{1.8} is finitely generated within the ring $\Sym\!\Bbb Q[M,L]$.
Eight polynomials \mythetag{1.9}, \mythetag{1.10}, \mythetag{1.11}, 
\mythetag{1.12}, \mythetag{1.13}, \mythetag{1.14}, \mythetag{1.15}, 
\mythetag{1.16} belong to the ideal $I_{\text{PC\kern -0.7pt\_\kern 0.5pt sym}}$ 
and constitute a basis of this ideal.
\endproclaim
    The theorem~\mythetheorem{1.1} was proved in \mycite{68}. Using the
polynomials \mythetag{1.9}, \mythetag{1.10}, \mythetag{1.11}, 
\mythetag{1.12}, \mythetag{1.13}, \mythetag{1.14}, \mythetag{1.15}, 
\mythetag{1.16}, in \mycite{67} the following equations 
were written:
$$
\xalignat 4
&\hskip -2em
\tilde p_{\kern 1pt 1}=0,
&&\tilde p_{\kern 1pt 2}=0,&&\tilde p_{\kern 1pt 3}=0,&&\tilde p_{\kern 1pt 4}=0,
\qquad\\
\vspace{-1.7ex}
\mytag{1.17}\\
\vspace{-1.7ex}
&\hskip -2em
\tilde p_{\kern 1pt 5}=0,&&\tilde p_{\kern 1pt 6}=0,&&\tilde p_{\kern 1pt 7}=0,
&&\tilde p_{\kern 1pt 8}=0.\qquad
\endxalignat
$$
The equations \mythetag{1.17} are factor equations of the perfect cuboid equations 
\mythetag{1.1} in the sense of the definition~\mythedefinition{1.2}. The 
theorem~\mythetheorem{1.1} means that the equations \mythetag{1.17} constitute 
a complete system of the factor equations.\par
     The equations \mythetag{1.17} are derived from the equations \mythetag{1.1}.
Therefore each solution of the equations \mythetag{1.1} is a solution for
the factor equations \mythetag{1.17}. Generally speaking, the converse 
proposition is not valid. However, in \mycite{67} the following theorem was proved.
\mytheorem{1.2} Each integer or rational solution of the factor equations 
\mythetag{1.17} such that $x_1>0$, $x_2>0$, $x_3>0$, $d_1>0$, $d_2>0$, and 
$d_3>0$ is an integer or rational solution for the equations \mythetag{1.1}.
\endproclaim
     The theorem~\mythetheorem{1.2} means that the factor equations \mythetag{1.17} 
are equivalent to the original equations \mythetag{1.1} with respect to the main 
problem of finding perfect cuboids or proving their non-existence. In this paper
we continue studying the factor equations \mythetag{1.17} by transforming them
into so-called $E$-form. As a result below in section 4 we derive a single
biquadratic equation from the equations \mythetag{1.17}.\par
\head
$E$-form of multisymmetric polynomials. 
\endhead
     Multisymmetric polynomials from the ring $\Sym\!\Bbb Q[M,L]$ are similar 
to regular symmetric polynomials (see \mycite{69}). Like in the case of regular 
symmetric polynomials, there are elementary symmetric polynomials in 
$\Sym\!\Bbb Q[M,L]$:
$$
\gather
\hskip -2em
\aligned
&e_{\sssize [1,0]}=x_1+x_2+x_3,\\
&e_{\sssize [2,0]}=x_1\,x_2+x_2\,x_3+x_3\,x_1,\\
&e_{\sssize [3,0]}=x_1\,x_2\,x_3,
\endaligned
\kern 3em
\aligned
&e_{\sssize [0,1]}=d_1+d_2+d_3,\\
&e_{\sssize [0,2]}=d_1\,d_2+d_2\,d_3+d_3\,d_1,\\
&e_{\sssize [0,3]}=d_1\,d_2\,d_3,
\endaligned
\quad
\mytag{2.1}\\
\vspace{2ex}
\hskip -2em
\aligned
&e_{\sssize [2,1]}=x_1\,x_2\,d_3+x_2\,x_3\,d_1+x_3\,x_1\,d_2,\\
&e_{\sssize [1,1]}=x_1\,d_2+d_1\,x_2+x_2\,d_3+d_2\,x_3+x_3\,d_1+d_3\,x_1,\\
&e_{\sssize [1,2]}=x_1\,d_2\,d_3+x_2\,d_3\,d_1+x_3\,d_1\,d_2.
\endaligned
\mytag{2.2}
\endgather
$$
The role of the polynomials \mythetag{2.1} and \mythetag{2.2} is described
by the following theorem. 
\mytheorem{2.1} The elementary multisymmetric polynomials \mythetag{2.1} and
\mythetag{2.2} generate the ring $\Sym\!\Bbb Q[M,L]$, i\.\,e\. each multisymmetric 
polynomial $p\in\Sym\!\Bbb Q[M,L]$ can be expressed as a polynomial with rational 
coefficients through these elementary multisymmetric polynomials.
\endproclaim
     The theorem~\mythetheorem{2.1} is known as the fundamental theorem for 
elementary multisymmetric polynomials. Its proof can be found in \mycite{52}.
\par
     Let's denote through $\Bbb Q[E,L]$ the following polynomial ring of ten 
independent variables $E_{10}$, $E_{20}$, $E_{30}$, $E_{01}$, $E_{02}$, $E_{03}$, 
$E_{21}$, $E_{11}$, $E_{12}$, $L$:
$$
\hskip -2em
\Bbb Q[E,L]=\Bbb Q[E_{10},E_{20},E_{30},E_{01},E_{02},E_{03},E_{21},E_{11},
E_{12},L].
\mytag{2.3}
$$
Like \mythetag{1.6}, the notation \mythetag{2.3} is used for the sake of brevity.
In terms of the notation \mythetag{2.3} the theorem~\mythetheorem{2.1} can be
formulated as follows. 
\mytheorem{2.2} For each multisymmetric polynomial $p\in\Sym\!\Bbb Q[M,L]$ there
is some polynomial $q\in\Bbb Q[E,L]$ such that $p$ is produced from $q$ by
substituting $e_{\sssize [1,0]}$, $e_{\sssize [2,0]}$, $e_{\sssize [3,0]}$,
$e_{\sssize [0,1]}$, $e_{\sssize [0,2]}$, $e_{\sssize [0,3]}$, $e_{\sssize [2,1]}$,
$e_{\sssize [1,1]}$, $e_{\sssize [1,2]}$ for $E_{10}$, $E_{20}$, $E_{30}$, $E_{01}$, 
$E_{02}$, $E_{03}$, $E_{21}$, $E_{11}$, $E_{12}$ into the arguments of the polynomial
$q$:
$$
p=q(e_{\sssize [1,0]},e_{\sssize [2,0]},e_{\sssize [3,0]},e_{\sssize [0,1]},
e_{\sssize [0,2]},e_{\sssize [0,3]},e_{\sssize [2,1]},e_{\sssize [1,1]},
e_{\sssize [1,2]},L).
\mytag{2.4}
$$
\endproclaim
     The substitution procedure \mythetag{2.4} determines a mapping:
$$
\hskip -2em
\varphi\!:\,\Bbb Q[E,L]\longrightarrow\Sym\!\Bbb Q[M,L].
\mytag{2.5}
$$ 
It is easy to see that the mapping \mythetag{2.5} is a ring homomorphism. Such a  
homomorphism is called a {\it substitution homomorphism}. The 
theorem~\mythetheorem{2.2} means that the substitution homomorphism \mythetag{2.5}
is surjective. 
\mydefinition{2.1} A polynomial $q=q(E_{10},E_{20},E_{30},E_{01},E_{02},E_{03},
E_{21},E_{11},E_{12},L)$ such that $p=\varphi(q)$ is called an $E$-form of a
polynomial $p\in\Sym\!\Bbb Q[M,L]$.
\enddefinition
     Unfortunately the substitution homomorphism \mythetag{2.5} is not bijective.
Therefore an $E$-form of a polynomial $p\in\Sym\!\Bbb Q[M,L]$ is not unique. It is
defined up to a polynomial from the kernel of the homomorphism \mythetag{2.5}. The
kernel $K=\Ker\varphi$ was studied in \mycite{68}. There it was shown that the 
ideal $K=\Ker\varphi$ has a Gr\"obner basis consisting of $14$ polynomials, i\.\,e\.
it is presented as 
$$
\hskip -2em
K=\bigl<\tilde q_{\kern 1pt 1},\tilde q_{\kern 1pt 2},\tilde q_{\kern 1pt 3},
\tilde q_{\kern 1pt 4},\tilde q_{\kern 1pt 5},\tilde q_{\kern 1pt 6},
\tilde q_{\kern 1pt 7},\tilde q_{\kern 1pt 8},\tilde q_{\kern 1pt 9},
\tilde q_{\kern 1pt 10},\tilde q_{\kern 1pt 11},\tilde q_{\kern 1pt 12},
\tilde q_{\kern 1pt 13},\tilde q_{\kern 1pt 14}\bigr>.
\mytag{2.6}
$$
For the definition of Gr\"obner bases and for their applications the reader is 
referred to the book \mycite{70}. The polynomials $\tilde q_{\kern 1pt 1}$, 
\dots, $\tilde q_{\kern 1pt 14}$ in \mythetag{2.6} were calculated with the use 
of symbolic computations, but they were not given in an explicit form in \mycite{68}.
Instead, another basis of the ideal $K$ consisting of seven polynomials was
suggested:
$$
\hskip -2em
K=\bigl<q_{\kern 1pt 1},q_{\kern 1pt 2},q_{\kern 1pt 3},q_{\kern 1pt 4},
q_{\kern 1pt 5},q_{\kern 1pt 6},q_{\kern 1pt 7}\bigr>.
\mytag{2.7}
$$
The explicit formulas for the polynomials $q_{\kern 1pt 1}$, \dots, 
$q_{\kern 1pt 7}$ from \mythetag{2.7} are available in \mycite{68}. 
The explicit expressions for $\tilde q_{\kern 1pt 1}$, \dots, $\tilde q_{\kern 1pt 14}$ 
are given in \mytheanchortext{App}{Appendix}.
\head
3. $E$-form of the factor equations. 
\endhead
     Now we are ready to proceed with studying the factor equations 
\mythetag{1.17}. For this purpose we transform them into the $E$-form as 
declared in the definition~\mythedefinition{2.1}. The first equation 
\mythetag{1.17} transformed into the $E$-form is very simple:
$$
\hskip -2em
E_{10}^2-2\,E_{20}-L^2=0.
\mytag{3.1}
$$
The second equation \mythetag{1.17} is a little bit more complicated:
$$
\hskip -2em
2\,E_{02}-4\,E_{20}-E_{01}^2+2\,E_{10}^2=0.
\mytag{3.2}
$$
Here are the third and the fourth equations \mythetag{1.17} transformed into
the $E$-form:
$$
\hskip -2em
\aligned
&E_{10}\,E_{11}-3\,E_{03}-E_{21}+3\,E_{01}\,E_{02}-E_{20}\,E_{01}-E_{01}^3=0,\\
&E_{01}\,E_{11}-E_{12}-3\,E_{30}+E_{10}\,E_{02}+E_{20}\,E_{10}-E_{01}^2\,E_{10}=0.
\endaligned
\mytag{3.3}
$$
Then we transform the fifth equation \mythetag{1.17}. As a result we get
$$
\hskip -2em
\aligned
-E_{10}\,&E_{21}-E_{01}\,E_{12}-E_{01}\,E_{30}-E_{01}^3\,E_{10}+E_{01}^2\,E_{11}\,-\\
&-\,E_{02}\,E_{11}+E_{11}\,E_{20}-E_{10}\,E_{03}+2\,E_{10}\,E_{01}\,E_{02}=0.
\endaligned
\mytag{3.4}
$$
The next step is to transform the sixth and the seventh equations \mythetag{1.17}.
Upon doing it we multiply both equations by $3$. Then we have 
$$
\gather
\hskip -2em
\gathered
4\,E_{01}\,E_{10}\,E_{11}-3\,E_{01}^2\,E_{10}^2+2\,E_{10}^2\,E_{02}
+2\,E_{20}\,E_{01}^2-2\,E_{10}\,E_{12}\,-\\
-\,2\,E_{02}\,E_{20}-2\,E_{01}\,E_{21}-E_{11}^2-12\,E_{10}\,E_{30}+6\,E_{20}^2=0.
\endgathered
\mytag{3.5}\\
\vspace{2ex}
\hskip -2em
\gathered
4\,E_{01}\,E_{10}\,E_{11}-4\,E_{10}^2\,E_{02}-4\,E_{20}\,E_{01}^2-2\,E_{10}\,E_{12}
+10\,E_{02}\,E_{20}\,-\\
-\,2\,E_{01}\,E_{21}-E_{11}^2-12\,E_{01}\,E_{03}-3\,E_{01}^4-6\,E_{02}^2
+12\,E_{01}^2\,E_{02}=0.
\endgathered
\mytag{3.6}
\endgather
$$
The last step is to transform the eighth equation \mythetag{1.17}. \pagebreak This 
equation is the most complicated of all eight. Upon transforming it we multiply 
this equation by $3$:
$$
\gathered
9\,E_{01}\,E_{03}\,E_{20}-7\,E_{01}^2\,E_{02}\,E_{20}+2\,E_{02}\,E_{10}\,E_{12}
-2\,E_{01}^2\,E_{10}\,E_{12}\,+\\
+\,3\,E_{03}\,E_{10}\,E_{11}+4\,E_{01}^3\,E_{10}\,E_{11}
-7\,E_{01}\,E_{02}\,E_{10}\,E_{11}-6\,E_{01}\,E_{03}\,E_{10}^2\,+\\
+\,8\,E_{01}^2\,E_{02}\,E_{10}^2+3\,E_{01}\,E_{11}\,E_{30}-2\,E_{01}\,E_{20}\,E_{21}
+E_{10}\,E_{12}\,E_{20}\,-\\
-\,E_{02}\,E_{10}^2\,E_{20}+E_{01}\,E_{10}\,E_{11}\,E_{20}+9\,E_{02}\,E_{10}\,E_{30}
-2\,E_{02}\,E_{20}^2\,+\\
+\,2\,E_{01}^2\,E_{20}^2-E_{11}^2\,E_{20}-3\,E_{12}\,E_{30}+E_{02}\,E_{11}^2
-E_{01}^2\,E_{11}^2\,-\\
-\,2\,E_{02}^2\,E_{10}^2+2\,E_{01}^4\,E_{20}+2\,E_{02}^2\,E_{20}-3\,E_{03}\,E_{21}\,-\\
-\,2\,E_{01}^3\,E_{21}+5\,E_{01}\,E_{02}\,E_{21}-6\,E_{01}^2\,E_{10}\,E_{30}
-3\,E_{01}^4\,E_{10}^2=0.
\endgathered\quad
\mytag{3.7}
$$
The transformed factor equations \mythetag{3.1}, \mythetag{3.2}, \mythetag{3.3}, 
\mythetag{3.4}, \mythetag{3.5}, \mythetag{3.6}, \mythetag{3.7} should be complemented
with the following $14$ kernel equations:
$$
\xalignat 3
&\hskip -2em
\tilde q_{\kern 1pt 1}=0,
&&\tilde q_{\kern 1pt 2}=0,
&&\tilde q_{\kern 1pt 3}=0,\\
&\hskip -2em
\tilde q_{\kern 1pt 4}=0,
&&\tilde q_{\kern 1pt 5}=0,
&&\tilde q_{\kern 1pt 6}=0,\\
&\hskip -2em
\tilde q_{\kern 1pt 7}=0,
&&\tilde q_{\kern 1pt 8}=0,
&&\tilde q_{\kern 1pt 9}=0,
\mytag{3.8}\\
&\hskip -2em
\tilde q_{\kern 1pt 10}=0,
&&\tilde q_{\kern 1pt 11}=0,
&&\tilde q_{\kern 1pt 12}=0,\\
&\hskip -2em
\tilde q_{\kern 1pt 13}=0,
&&\tilde q_{\kern 1pt 14}=0.
\endxalignat
$$
The kernel polynomials $\tilde q_{\kern 1pt 1}$, \dots, $\tilde q_{\kern 1pt 14}$ 
are rather huge. The explicit expressions for them are given in 
\mytheanchortext{App}{Appendix} in a machine readable form.\par
     The equations \mythetag{3.1}, \mythetag{3.2}, \mythetag{3.3}, \mythetag{3.4}, 
\mythetag{3.5}, \mythetag{3.6}, \mythetag{3.7} complemented with the equations
\mythetag{3.8} constitute a huge system of $22$ polynomial equations with integer
coefficients with respect to $10$ variables. Despite being a system of polynomial 
equations with integer coefficients, we can consider real solutions of this 
system. These real solutions constitute a real algebraic variety in $\Bbb R^{10}$. 
We denote this variety through $E^{10}_{\text{PC\kern -0.7pt\_\kern 0.5pt sym}}$. 
Integer points of this variety, provided they do exist, are associated with 
perfect cuboids. The real algebraic variety $E^{10}_{\text{PC\kern -0.7pt\_\kern 0.5pt 
sym}}\subset\Bbb R^{10}$ is studied below in section 4.\par
\head
4. Reduction of the factor equations. 
\endhead
     Let's consider the equation \mythetag{3.1}. This equation is linear with
respect to the variable $E_{20}$. Resolving it with respect to this variable,
we get
$$
\hskip -2em
E_{20}=\frac{1}{2}\,E_{10}^2-\frac{1}{2}\,L^2.
\mytag{4.1}
$$
Substituting \mythetag{4.1} into \mythetag{3.2}, we derive 
$$
\hskip -2em
2\,L^2+2\,E_{02}-E_{01}^2=0.
\mytag{4.2}
$$
The equation \mythetag{4.2} is similar to \mythetag{3.1}. This equation is linear 
with respect to the variable $E_{02}$. Resolving it with respect to this variable,
we get
$$
\hskip -2em
E_{02}=\frac{1}{2}\,E_{01}^2-L^2.
\mytag{4.3}
$$
Now we substitute \mythetag{4.1} and \mythetag{4.3} into \mythetag{3.3}. This yields 
the equations 
$$
\align
&\hskip -2em
-3\,E_{03}-E_{21}-\frac{5}{2}\,E_{01}\,L^2+\frac{1}{2}\,E_{01}^3+E_{10}\,E_{11}
-\frac{1}{2}\,E_{01}\,E_{10}^2=0,\quad
\mytag{4.4}\\
\vspace{1ex}
&\hskip -2em
-3\,E_{30}-E_{12}-\frac{3}{2}\,E_{10}\,L^2-\frac{1}{2}\,E_{01}^2\,E_{10}
+E_{01}\,E_{11}+\frac{1}{2}\,E_{10}^3=0.\quad
\mytag{4.5}
\endalign
$$
The equation \mythetag{4.4} is linear with respect to $E_{03}$, while the equation
\mythetag{4.5} is linear with respect to $E_{30}$. Resolving these equations, we derive
$$
\align
&\hskip -2em
E_{03}=-\frac{1}{3}\,E_{21}-\frac{1}{6}\,E_{01}\,E_{10}^2-\frac{5}{6}\,E_{01}\,L^2
+\frac{1}{6}\,E_{01}^3+\frac{1}{3}\,E_{10}\,E_{11},
\mytag{4.6}\\
\vspace{1ex}
&\hskip -2em
E_{30}=-\frac{1}{3}\,E_{12}-\frac{1}{6}\,E_{10}\,E_{01}^2
-\frac{1}{2}\,E_{10}\,L^2+\frac{1}{6}\,E_{10}^3
+\frac{1}{3}\,E_{01}\,E_{11}.
\mytag{4.7}
\endalign
$$
Note that the formulas \mythetag{4.1}, \mythetag{4.3}, \mythetag{4.6}, and
\mythetag{4.7} were already derived in \mycite{45}.\par
     Using the formulas \mythetag{4.1}, \mythetag{4.3}, \mythetag{4.6}, \mythetag{4.7}, 
we can eliminate the variables $E_{20}$, $E_{02}$, $E_{30}$, $E_{03}$ from the remaining
factor equations \mythetag{3.4}, \mythetag{3.5}, \mythetag{3.6}, \mythetag{3.7} and
from the kernel equations \mythetag{3.8}. As a result the number of variables reduces
from $10$ to $6$, while the total number of equations reduced from $22$ to $18$. The 
reduced system of $18$ equations determines a real algebraic variety in $\Bbb R^6$. We
denote this variety through $E^6_{\text{PC\kern -0.7pt\_\kern 0.5pt sym}}$. The algebraic
variety $E^6_{\text{PC\kern -0.7pt\_\kern 0.5pt sym}}\subset\Bbb R^6$ is the projection
of the variety $E^{10}_{\text{PC\kern -0.7pt\_\kern 0.5pt sym}}\subset\Bbb R^{10}$ onto
the subspace $\Bbb R^6\subset\Bbb R^{10}$:
$$
\hskip -2em
\pi_{6}\!:\,E^{10}_{\text{PC\kern -0.7pt\_\kern 0.5pt sym}}\to
E^6_{\text{PC\kern -0.7pt\_\kern 0.5pt sym}}.
\mytag{4.8}
$$ 
The right hand sides of the formulas \mythetag{4.1}, \mythetag{4.3}, \mythetag{4.6}, 
\mythetag{4.7} are polynomials with respect to six variables $E_{10}$, $E_{01}$, 
$E_{21}$, $E_{11}$, $E_{12}$, and $L$. For this reason the mapping \mythetag{4.8} is
bijective, i\.\,e\. there is the inverse mapping 
$$
\hskip -2em
\pi_{6}^{-1}\!:\,E^6_{\text{PC\kern -0.7pt\_\kern 0.5pt sym}}\to
E^{10}_{\text{PC\kern -0.7pt\_\kern 0.5pt sym}}.
\mytag{4.9}
$$\par
     As it was said above, the equations \mythetag{3.4}, \mythetag{3.5}, \mythetag{3.6}, 
\mythetag{3.7} are transformed with the use of \mythetag{4.1}, \mythetag{4.3}, \mythetag{4.6},
\mythetag{4.7}. Here is the explicit form of the transformed equation \mythetag{3.4}
multiplied by the number $-3/2$:
$$
E_{10}\,E_{21}+E_{01}\,E_{12}=\frac{1}{4}\,E_{11}\,E_{10}^2
+\frac{1}{4}\,E_{01}^2\,E_{11}+\frac{3}{4}\,E_{11}\,L^2-E_{10}\,E_{01}\,L^2.\quad
\mytag{4.10}
$$
Two terms of the equation \mythetag{4.10} are left in its left hand side, the other four
terms are brought to the right hand side of this equation.\par
     Upon transforming the equations \mythetag{3.5} and \mythetag{3.6} with the use of
\mythetag{4.1}, \mythetag{4.3}, \mythetag{4.6}, \mythetag{4.7} we do not write them separately,
but compose their sum and their difference. The difference of the equations \mythetag{3.5} and 
\mythetag{3.6} transformed by means of \mythetag{4.1}, \mythetag{4.3}, \mythetag{4.6}, 
\mythetag{4.7} and then multiplied by $1/4$ looks like 
$$
-E_{01}\,E_{21}+E_{10}\,E_{12}=\frac{1}{8}\,E_{10}^4-\frac{1}{8}\,E_{01}^4
-\frac{3}{4}\,E_{10}^2\,L^2+E_{01}^2\,L^2-\frac{3}{8}\,L^4.\quad
\mytag{4.11}
$$
Like in \mythetag{4.10}, Two terms of the equation \mythetag{4.11} are left in its left 
hand side, the other five terms are moved to the right hand side of this equation.\par
     The sum of the equations \mythetag{3.5} and \mythetag{3.6} transformed by means of 
\mythetag{4.1}, \mythetag{4.3}, \mythetag{4.6}, \mythetag{4.7} and then multiplied by 
the number $-2$ is written as
$$
4\,E_{11}^2+E_{10}^4+E_{01}^4-2\,E_{10}^2\,E_{01}^2-2\,L^2\,E_{10}^2
-6\,E_{01}^2\,L^2+L^4=0.\quad
\mytag{4.12}
$$
And finally, the equation \mythetag{3.7} transformed by means of 
\mythetag{4.1}, \mythetag{4.3}, \mythetag{4.6}, \mythetag{4.7} and then multiplied
by the number $4$ is written as
$$
\gathered
8\,E_{10}\,E_{11}\,E_{21}+8\,E_{01}\,E_{11}\,E_{12}-4\,E_{21}^2-4\,E_{12}^2
-8\,E_{10}\,E_{12}\,L^2=\\
=2\,E_{10}^2\,E_{11}^2+2\,E_{01}^2\,E_{10}^2\,L^2
+2\,E_{01}^2\,E_{11}^2-E_{01}^4\,L^2-2\,E_{10}^4\,L^2\,-\\
-\,8\,E_{01}\,E_{10}\,E_{11}\,L^2+8\,E_{10}^2\,L^4+6\,E_{01}^2\,L^4
-2\,E_{11}^2\,L^2-2\,L^6.
\endgathered
\quad
\mytag{4.13}
$$
It is easy to see that the transformed equation \mythetag{4.13} looks much more 
simple than the original equation \mythetag{3.7}.\par
     The kernel equations \mythetag{3.8} also get more simple when transformed by
means of \mythetag{4.1}, \mythetag{4.3}, \mythetag{4.6}, \mythetag{4.7}. But they 
are still rather huge for to write them explicitly.\par
\head
5. Further transformations. 
\endhead
     Note that the equations \mythetag{4.10} and \mythetag{4.11} are linear
with respect to the variables $E_{21}$ and $E_{12}$. They can be resolved with
respect to these variables as
$$
\align
&\hskip -2em
\aligned
E_{21}&=\frac{2\,E_{10}^3\,E_{11}+2\,E_{01}^2\,E_{10}\,E_{11}
-E_{01}\,E_{10}^4+E_{01}^5}
{8\,(E_{01}^2+E_{10}^2)\vphantom{\vrule height 11pt}}\,+\\
\vspace{1ex}
&+\,\frac{6\,E_{10}\,E_{11}\,L^2-2\,E_{01}\,E_{10}^2\,L^2
-8\,E_{01}^3\,L^2+3\,E_{01}\,L^4}
{8\,(E_{01}^2+E_{10}^2)\vphantom{\vrule height 11pt}},
\endaligned
\mytag{5.1}\\
\vspace{2ex}
&\hskip -2em
\aligned
E_{12}=&\frac{E_{01}^4\,E_{10}-2\,E_{01}^3\,E_{11}
-2\,E_{01}\,E_{10}^2\,E_{11}-E_{10}^5}
{8\,(E_{01}^2+E_{10}^2)\vphantom{\vrule height 11pt}}\,+\\
\vspace{1ex}
&\kern 5em +\,\frac{6\,E_{10}^3\,L^2-6\,E_{01}\,E_{11}\,L^2+3\,E_{10}\,L^4}
{8\,(E_{01}^2+E_{10}^2)\vphantom{\vrule height 11pt}}.
\endaligned
\mytag{5.2}
\endalign
$$
The equation \mythetag{4.12} can be resolved with respect to the square of
$E_{11}$:
$$
\hskip -2em
E_{11}^2=\frac{1}{2}\,E_{01}^2\,E_{10}^2-\frac{1}{4}\,E_{10}^4
-\frac{1}{4}\,E_{01}^4+\frac{1}{2}\,E_{10}^2\,L^2+\frac{3}{2}\,E_{01}^2\,L^2
-\frac{1}{4}\,L^4. 
\mytag{5.3}
$$\par
     Using \mythetag{5.1} and \mythetag{5.2}, we can eliminate $E_{21}$ and
$E_{12}$ from the remaining factor equation \mythetag{4.13} and from the kernel
equations \mythetag{3.8} which are already transformed with the use of 
\mythetag{4.1}, \mythetag{4.3}, \mythetag{4.6}, \mythetag{4.7}. Then we can apply 
\mythetag{5.3} and eliminate squares of $E_{11}$, cubes of $E_{11}$, and all other
higher degrees of $E_{11}$ from \mythetag{4.13} and \mythetag{3.8}. As a result of 
such transformations the equations \mythetag{4.13} and \mythetag{3.8} luckily
turn to trivial identities $0=0$.\par
     Note that the formulas \mythetag{5.1} and \mythetag{5.2} have nontrivial
denominators. Therefore  we should be careful in interpreting our lucky result. 
The determinants in \mythetag{5.1} and \mythetag{5.2} turn to zero at the points
where $E_{10}=0$ and $E_{01}=0$ simultaneously. Hence we can formulate the 
following theorem.
\mytheorem{5.1} At the locus of points in $\Bbb R^{10}$ where $E_{10}\neq 0$
or  $E_{01}\neq 0$ the system of\/ $22$ polynomial equations \mythetag{3.1}, 
\mythetag{3.2}, \mythetag{3.3}, \mythetag{3.4}, \mythetag{3.5}, \mythetag{3.6}, 
\mythetag{3.7}, and \mythetag{3.8} is equivalent to the single equation
\mythetag{4.12}.
\endproclaim
     The equation \mythetag{4.12} is a polynomial equation with respect to 
four variables $E_{10}$, $E_{01}$, $E_{11}$, and $L$. It defines a real algebraic
variety in $\Bbb R^4$. Let's denote it through $E^4_{\text{PC\kern -0.7pt\_\kern 0.5pt 
sym}}$. Like $E^{10}_{\text{PC\kern -0.7pt\_\kern 0.5pt sym}}$ and $E^6_{\text{PC\kern 
-0.7pt\_\kern 0.5pt sym}}$, the variety $E^4_{\text{PC\kern -0.7pt\_\kern 0.5pt sym}}$
can have special points where $E_{10}=0$ and $E_{01}=0$ simultaneously. Let's pin out
these special points from each of the three varieties. As a result we get three
Zariski open subsets $\tilde E^{10}_{\text{PC\kern -0.7pt\_\kern 0.5pt sym}}$, 
$\tilde E^6_{\text{PC\kern -0.7pt\_\kern 0.5pt sym}}$, and $\tilde E^4_{\text{PC\kern 
-0.7pt\_\kern 0.5pt sym}}$ within $E^{10}_{\text{PC\kern -0.7pt\_\kern 0.5pt sym}}$, 
$E^6_{\text{PC\kern -0.7pt\_\kern 0.5pt sym}}$, and $E^4_{\text{PC\kern -0.7pt\_\kern 
0.5pt sym}}$ respectively.\par
     Generally speaking, the algebraic variety $E^4_{\text{PC\kern -0.7pt\_\kern 
0.5pt sym}}$ should not be a projection of $E^6_{\text{PC\kern -0.7pt\_\kern 0.5pt sym}}$.
However, its subset $\tilde E^4_{\text{PC\kern -0.7pt\_\kern 0.5pt sym}}$ is the
projection of the subset $\tilde E^6_{\text{PC\kern -0.7pt\_\kern 0.5pt sym}}$ from
$\Bbb R^6$ onto the subspace $\Bbb R^4\subset\Bbb R^6$. We have the mapping
$$
\hskip -2em
\pi_{4}\!:\,\tilde E^6_{\text{PC\kern -0.7pt\_\kern 0.5pt sym}}\to
\tilde E^4_{\text{PC\kern -0.7pt\_\kern 0.5pt sym}}.
\mytag{5.4}
$$ 
Due to \mythetag{5.1} and \mythetag{5.2} the mapping \mythetag{5.4} is bijective, 
i\.\,e\. it has the inverse mapping 
$$
\hskip -2em
\pi_{4}^{-1}\!:\,\tilde E^4_{\text{PC\kern -0.7pt\_\kern 0.5pt sym}}\to
\tilde E^6_{\text{PC\kern -0.7pt\_\kern 0.5pt sym}}.
\mytag{5.5}
$$
The mappings \mythetag{5.4} and \mythetag{5.5} are similar to the mappings
\mythetag{4.8} and \mythetag{4.9}.\par
     The equation \mythetag{4.12} was derived from the equations \mythetag{3.1}, 
\mythetag{3.2}, \mythetag{3.3}, \mythetag{3.4}, \mythetag{3.5}, \mythetag{3.6}, 
\mythetag{3.7}, and \mythetag{3.8} without use of the denominators in
\mythetag{5.1} and \mythetag{5.2}. For this reason each integer solution of 
the equations \mythetag{3.1}, \mythetag{3.2}, \mythetag{3.3}, \mythetag{3.4}, 
\mythetag{3.5}, \mythetag{3.6}, \mythetag{3.7}, \mythetag{3.8} induces an
integer solution of the equation \mythetag{4.12}. The backward formulas 
\mythetag{4.1}, \mythetag{4.3}, \mythetag{4.6}, \mythetag{4.7}, \mythetag{5.1}, 
\mythetag{5.2} comprise fractions. Nevertheless, we can prove the following theorem.
\mytheorem{5.2} Each integer or rational solution of the equation \mythetag{4.12}
such that $E_{10}^2+E_{01}^2\neq 0$ induces an integer solution of the complete 
system of $22$ equations \mythetag{3.1}, \mythetag{3.2}, \mythetag{3.3}, 
\mythetag{3.4}, \mythetag{3.5}, \mythetag{3.6}, \mythetag{3.7}, \mythetag{3.8}.
\endproclaim
\demo{Proof} The equations \mythetag{3.1}, \mythetag{3.2}, \mythetag{3.3}, 
\mythetag{3.4}, \mythetag{3.5}, \mythetag{3.6}, \mythetag{3.7}, \mythetag{3.8}
are weighted homogeneous. If we have a solution of them, then applying the
transformations 
$$
\xalignat 2
&\hskip -2em
E_{10}\mapsto\alpha\,E_{10},
&&E_{01}\mapsto\alpha\,E_{01},\\
&\hskip -2em
E_{20}\mapsto\alpha^2\,E_{20},
&&E_{02}\mapsto\alpha^2\,E_{02},\\
&\hskip -2em
E_{30}\mapsto\alpha^3\,E_{30},
&&E_{03}\mapsto\alpha^3\,E_{03},
\mytag{5.6}\\
&\hskip -2em
E_{21}\mapsto\alpha^3\,E_{21},
&&E_{12}\mapsto\alpha^3\,E_{12},\\
&\hskip -2em
E_{11}\mapsto\alpha^2\,E_{11},
&&L\mapsto\alpha\,L,\\
\endxalignat
$$
we get another solution for them. Assume that we have an integer or a rational
solution of the equation \mythetag{4.12} such that $E_{10}^2+E_{01}^2\neq 0$. 
Then by means of the formulas \mythetag{4.1}, \mythetag{4.3}, \mythetag{4.6}, 
\mythetag{4.7}, \mythetag{5.1}, \mythetag{5.2} we get a solution for the equations
\mythetag{3.1}, \mythetag{3.2}, \mythetag{3.3}, \mythetag{3.4}, \mythetag{3.5}, 
\mythetag{3.6}, \mythetag{3.7}, \mythetag{3.8} such that $E_{10}$, $E_{20}$,
$E_{30}$, $E_{01}$, $E_{02}$, $E_{03}$, $E_{21}$, $E_{11}$, $E_{12}$, and $L$
are integer or rational numbers. Let $Q$ be a common denominator for all of these
ten numbers. It is sufficient to apply the transformations \mythetag{5.6}
with $\alpha=Q$ and obtain an integer solution for the equations \mythetag{3.1}, 
\mythetag{3.2}, \mythetag{3.3}, \mythetag{3.4}, \mythetag{3.5}, \mythetag{3.6}, 
\mythetag{3.7}, \mythetag{3.8}.
\qed\enddemo
\head
6. Back to perfect cuboids. 
\endhead
     If a perfect Euler cuboid does exist, it induces an integer solution for 
the equations \mythetag{1.1} and an integer solution for the factor equations 
\mythetag{1.17} such that $x_1>0$, $x_2>0$, $x_3>0$, $d_1>0$, $d_2>0$, $d_1>0$, 
and $L>0$. Via the formulas \mythetag{2.1} and \mythetag{2.2} the latter one 
induces a positive integer solution for the equations \mythetag{3.1}, \mythetag{3.2}, 
\mythetag{3.3}, \mythetag{3.4}, \mythetag{3.5}, \mythetag{3.6}, \mythetag{3.7}, 
\mythetag{3.8}, and finally, an integer solution for the equation \mythetag{4.12}
with $E_{10}>0$, $E_{01}>0$, $E_{11}>0$ and $L>0$. Combining this fact with the 
theorem~\mythetheorem{5.2} we formulate the following result.
\mytheorem{6.1} The existence of positive rational solutions of the equation 
\mythetag{4.12} is a necessary condition for the existence of perfect Euler
cuboids. 
\endproclaim
     The backward path from the equation \mythetag{4.12} to perfect cuboids is
not so straightforward. On this path we need to solve the inverse problem of 
calculating the values of $x_1$, $x_2$, $x_3$, $d_1$, $d_2$, $d_3$ through the 
known values of the elementary multisymmetric polynomials in \mythetag{2.1} and 
\mythetag{2.2}. As I know, this problem is poorly studied. Therefore the 
theorem~\mythetheorem{6.1} is formulated as a necessary, but not a sufficient,
condition.\par 
     The theorem~\mythetheorem{6.1} is the main result of this paper. It could
be useful in computerized search for perfect cuboids or in proving their 
non-existence.\par
\head
7. Similarity to Heron's problem. 
\endhead
     The ancient Greek mathematician Heron of Alexandria (10--70 C\.E\.) has 
discovered the formula expressing the area $S$ of a triangle through its sides 
$a$, $b$, and $c$. This Heron's formula is written as follows:
$$
\hskip -2em
S=\sqrt{p\,(p-a)\,(p-b)\,(p-c)}\text{, \ where \ }p=\frac{a+b+c}{2}.
\mytag{7.1}
$$
The formula \mythetag{7.1} can be written as an equation with respect to $a$, 
$b$, $c$, and $S$:
$$
\hskip -2em
(4\,S)^2+(a^2+b^2-c^2)^2-4\,a^2\,b^2=0.
\mytag{7.2}
$$     
The equation \mythetag{7.2} is associated with the well known Heron's problem
--- find all triangles whose sides and whose area are integer numbers (see
\mycite{71}).\par
     Now let's consider the equation \mythetag{4.12}. This equation can be 
rewritten as 
$$
\hskip -2em
(2\,E_{11})^2+(E_{01}^2+L^2-E_{10}^2)^2-8\,E_{01}^2\,L^2=0. 
\mytag{7.3}
$$
The equation \mythetag{7.3} is very similar to the equation \mythetag{7.2}.
Both equations belong to the class of biquadratic Diophantine equations. 

\par
\head
\myanchortext{App}{Appendix.}
\endhead
Here are the formulas for the polynomials $\tilde q_{\kern 1pt 1}$, \dots, 
$\tilde q_{\kern 1pt 14}$ written in a machine readable form convenient for to
copy-paste into some symbolic computations package:
\medskip\noindent
{\tt\~q1:=-3*E02*E21+E01\^{}2*E21-9*E03*E20+4*E01*E02*E20\relax
-E01\^{}3*E20+3*E11*E12\newline
-2*E01*E10*E12-E01*E11\^{}2+E01\^{}2*E10*E11+3*E03*E10\^{}2-E01*E02*E10\^{}2;}
\medskip\noindent
{\tt\~q2:=9*E02*E30-3*E01\^{}2*E30-3*E11*E21+2*E01*E10*E21+3*E12*E20-4*E02\newline
*E10*E20+E01\^{}2*E10*E20-E10\^{}2*E12+E10*E11\^{}2-E01*E10\^{}2*E11+E02*E10\^{}3;}
\medskip\noindent
{\tt\~q3:=-9*E02\^{}2*E30+6*E01\^{}2*E02*E30-E01\^{}4*E30-3*E02*E12*E20+E01\^{}2%
*E12\newline
*E20-9*E03*E11*E20+4*E01*E02*E11*E20-E01\^{}3*E11*E20+6*E01*E03*E10*E20\newline
+4*E02\^{}2*E10*E20-5*E01\^{}2*E02*E10*E20+E01\^{}4*E10*E20+3*E11\^{}2*E12\relax
-4*E01\newline
*E10*E11*E12+E02*E10\^{}2*E12+E01\^{}2*E10\^{}2*E12-E01*E11\^{}3\relax
-E02*E10*E11\^{}2\newline
+2*E01\^{}2*E10*E11\^{}2+3*E03*E10\^{}2*E11-E01\^{}3*E10\^{}2*E11\relax
-2*E01*E03*E10\^{}3\newline
-E02\^{}2*E10\^{}3+E01\^{}2*E02*E10\^{}3;}
\medskip\noindent
{\tt\~q4:=-27*E03*E21+E01\^{}3*E21-18*E01*E03*E20+12*E02\^{}2*E20\relax
+E01\^{}2*E02*E20\newline
-E01\^{}4*E20+9*E12\^{}2+3*E01*E11*E12-6*E02*E10*E12-2*E01\^{}2*E10*E12\relax
-3*E02\newline
*E11\^{}2-E01\^{}2*E11\^{}2+9*E03*E10*E11+3*E01*E02*E10*E11\relax
+E01\^{}3*E10*E11+3\newline
*E01*E03*E10\^{}2-3*E02\^{}2*E10\^{}2-E01\^{}2*E02*E10\^{}2;}
\medskip\noindent
{\tt\~q5:=-81*E03*E30+18*E01*E02*E30-3*E01\^{}3*E30+9*E12*E21\relax
-E01\^{}2*E10*E21\newline
-6*E01*E12*E20+12*E02*E11*E20-3*E01\^{}2*E11*E20+36*E03*E10*E20-16*E01\newline
*E02*E10*E20+4*E01\^{}3*E10*E20-6*E10*E11*E12+5*E01*E10\^{}2*E12-3*E11\^{}3\newline
+7*E01*E10*E11\^{}2-3*E02*E10\^{}2*E11-4*E01\^{}2*E10\^{}2*E11\relax
-9*E03*E10\^{}3+4\newline
*E01*E02*E10\^{}3;}
\medskip\noindent
{\tt\~q6:=-9*E12*E30+3*E01*E11*E30-3*E02*E10*E30+3*E21\^{}2-2*E01*E20*E21+4\newline
*E02*E20\^{}2-E01\^{}2*E20\^{}2+E10*E12*E20-E11\^{}2*E20+E01*E10*E11*E20\relax
-E02\newline *E10\^{}2*E20;}
\medskip\noindent
{\tt\~q7:=-27*E02*E03*E30+9*E01\^{}2*E03*E30+3*E01*E02\^{}2*E30-E01\^{}\relax
3*E02*E30-9\newline
*E03*E12*E20+E01*E02*E12*E20-3*E01*E03*E11*E20+4*E02\^{}2*E11*E20-E01\^{}2\newline
*E02*E11*E20+12*E02*E03*E10*E20-E01\^{}2*E03*E10*E20-4*E01*E02\^{}2*E10\newline
*E20+E01\^{}3*E02*E10*E20+3*E11*E12\^{}2-2*E01*E10*E12\^{}2-2*E02*E10*E11*E12\newline
+3*E03*E10\^{}2*E12+E01*E02*E10\^{}2*E12-E02*E11\^{}3+2*E01*E02*E10*E11\^{}2+E01\newline
*E03*E10\^{}2*E11-E02\^{}2*E10\^{}2*E11-E01\^{}2*E02*E10\^{}2*E11-3*E02*E03*E10\^{}3\newline
+E01*E02\^{}2*E10\^{}3;}
\medskip\noindent
{\tt\~q8:=27*E02*E12*E30-9*E01\^{}2*E12*E30-81*E03*E11*E30+18*E01*E02*E11*E30\newline
-3*E01\^{}3*E11*E30+54*E01*E03*E10*E30-9*E02\^{}2*E10*E30-9*E01\^{}2*E02*E10\newline
*E30+2*E01\^{}4*E10*E30+9*E12\^{}2*E20-6*E01*E11*E12*E20-15*E02*E10*E12\newline
*E20+7*E01\^{}2*E10*E12*E20+12*E02*E11\^{}2*E20-3*E01\^{}2*E11\^{}2*E20+27*E03\newline
*E10*E11*E20-20*E01*E02*E10*E11*E20+5*E01\^{}3*E10*E11*E20-18*E01*E03\newline
*E10\^{}2*E20+4*E02\^{}2*E10\^{}2*E20+7*E01\^{}2*E02*E10\^{}2*E20\relax
-2*E01\^{}4*E10\^{}2*E20\newline
-3*E10\^{}2*E12\^{}2+2*E01*E10\^{}2*E11*E12+4*E02*E10\^{}3*E12\relax
-2*E01\^{}2*E10\^{}3*E12\newline
-3*E11\^{}4+8*E01*E10*E11\^{}3-4*E02*E10\^{}2*E11\^{}2-7*E01\^{}2*E10\^{}2*E11\^{}2\relax
-6*E03\newline
*E10\^{}3*E11+6*E01*E02*E10\^{}3*E11+2*E01\^{}3*E10\^{}3*E11+4*E01*E03*E10\^{}4\newline
-E02\^{}2*E10\^{}4-2*E01\^{}2*E02*E10\^{}4;}
\medskip\noindent
{\tt\~q9:=-27*E03\^{}2*E20+18*E01*E02*E03*E20-4*E01\^{}3*E03*E20-4*E02\^{}3*E20\newline
+E01\^{}2*E02\^{}2*E20-3*E02*E12\^{}2+E01\^{}2*E12\^{}2+9*E03*E11*E12-E01*E02*E11\newline
*E12-6*E01*E03*E10*E12+2*E02\^{}2*E10*E12-3*E01*E03*E11\^{}2+E02\^{}2*E11\^{}2\newline
-3*E02*E03*E10*E11+4*E01\^{}2*E03*E10*E11-E01*E02\^{}2*E10*E11+9*E03\^{}2\newline
*E10\^{}2-4*E01*E02*E03*E10\^{}2+E02\^{}3*E10\^{}2;}
\medskip\noindent
{\tt\~q10:=-243*E02*E03*E20*E30+81*E01\^{}2*E03*E20*E30+81*E01*E02\^{}2*E20*E30\newline
-45*E01\^{}3*E02*E20*E30+6*E01\^{}5*E20*E30-81*E02*E11*E12*E30+27*E01\^{}2\newline
*E11*E12*E30+54*E01*E02*E10*E12*E30-18*E01\^{}3*E10*E12*E30+243*E03\newline
*E11\^{}2*E30-54*E01*E02*E11\^{}2*E30+9*E01\^{}3*E11\^{}2*E30-324*E01*E03*E10\newline
*E11*E30+27*E02\^{}2*E10*E11*E30+63*E01\^{}2*E02*E10*E11*E30-12*E01\^{}4*E10\newline
*E11*E30+81*E02*E03*E10\^{}2*E30+81*E01\^{}2*E03*E10\^{}2*E30-45*E01*E02\^{}2\newline
*E10\^{}2*E30-3*E01\^{}3*E02*E10\^{}2*E30+2*E01\^{}5*E10\^{}2*E30-81*E03*E12*E20\^{}2\newline
+27*E01*E02*E12*E20\^{}2-6*E01\^{}3*E12*E20\^{}2+27*E01*E03*E11*E20\^{}2+36\newline
*E02\^{}2*E11*E20\^{}2-33*E01\^{}2*E02*E11*E20\^{}2+6*E01\^{}4*E11*E20\^{}2+108*E02\newline
*E03*E10*E20\^{}2-45*E01\^{}2*E03*E10*E20\^{}2-60*E01*E02\^{}2*E10*E20\^{}2+39\newline
*E01\^{}3*E02*E10*E20\^{}2-6*E01\^{}5*E10*E20\^{}2+27*E02*E10*E11*E12*E20-9\newline
*E01\^{}2*E10*E11*E12*E20+54*E03*E10\^{}2*E12*E20-36*E01*E02*E10\^{}2*E12*E20\newline
+10*E01\^{}3*E10\^{}2*E12*E20-45*E02*E11\^{}3*E20+15*E01\^{}2*E11\^{}3*E20-81*E03\newline
*E10*E11\^{}2*E20+108*E01*E02*E10*E11\^{}2*E20-33*E01\^{}3*E10*E11\^{}2*E20+90\newline
*E01*E03*E10\^{}2*E11*E20-33*E02\^{}2*E10\^{}2*E11*E20-59*E01\^{}2*E02*E10\^{}2*E11\newline
*E20+20*E01\^{}4*E10\^{}2*E11*E20-63*E02*E03*E10\^{}3*E20-9*E01\^{}2*E03*E10\^{}3\newline
*E20+43*E01*E02\^{}2*E10\^{}3*E20-5*E01\^{}3*E02*E10\^{}3*E20-2*E01\^{}5*E10\^{}3*E20\newline
-6*E02*E10\^{}3*E11*E12+2*E01\^{}2*E10\^{}3*E11*E12-9*E03*E10\^{}4*E12+7*E01*E02\newline
*E10\^{}4*E12-2*E01\^{}3*E10\^{}4*E12+9*E11\^{}5-30*E01*E10*E11\^{}4+15*E02*E10\^{}2\newline
*E11\^{}3+35*E01\^{}2*E10\^{}2*E11\^{}3+18*E03*E10\^{}3*E11\^{}2-34*E01*E02*E10\^{}3\newline
*E11\^{}2-16*E01\^{}3*E10\^{}3*E11\^{}2-21*E01*E03*E10\^{}4*E11+6*E02\^{}2*E10\^{}4*E11\newline
+21*E01\^{}2*E02*E10\^{}4*E11+2*E01\^{}4*E10\^{}4*E11+9*E02*E03*E10\^{}5+4*E01\^{}2\newline
*E03*E10\^{}5-7*E01*E02\^{}2*E10\^{}5-2*E01\^{}3*E02*E10\^{}5;}
\medskip\noindent
{\tt\~q11:=81*E02\^{}2*E12*E30-54*E01\^{}2*E02*E12*E30+9*E01\^{}4*E12*E30-243*E02\newline
*E03*E11*E30+81*E01\^{}2*E03*E11*E30+54*E01*E02\^{}2*E11*E30-27*E01\^{}3*E02\newline
*E11*E30+3*E01\^{}5*E11*E30+162*E01*E02*E03*E10*E30-54*E01\^{}3*E03*E10\newline
*E30-27*E02\^{}3*E10*E30-18*E01\^{}2*E02\^{}2*E10*E30+15*E01\^{}4*E02*E10*E30-2\newline
*E01\^{}6*E10*E30-243*E03\^{}2*E20\^{}2+162*E01*E02*E03*E20\^{}2-36*E01\^{}3*E03\newline
*E20\^{}2-36*E02\^{}3*E20\^{}2+9*E01\^{}2*E02\^{}2*E20\^{}2+81*E03*E11*E12*E20-27*E01\newline
*E02*E11*E12*E20+6*E01\^{}3*E11*E12*E20-54*E01*E03*E10*E12*E20-27*E02\^{}2\newline
*E10*E12*E20+36*E01\^{}2*E02*E10*E12*E20-7*E01\^{}4*E10*E12*E20-27*E01*E03\newline
*E11\^{}2*E20+45*E02\^{}2*E11\^{}2*E20-21*E01\^{}2*E02*E11\^{}2*E20+3*E01\^{}4*E11\^{}2\newline
*E20+54*E02*E03*E10*E11*E20+9*E01\^{}2*E03*E10*E11*E20-69*E01*E02\^{}2*E10\newline
*E11*E20+35*E01\^{}3*E02*E10*E11*E20-5*E01\^{}5*E10*E11*E20+162*E03\^{}2\newline
*E10\^{}2*E20-144*E01*E02*E03*E10\^{}2*E20+30*E01\^{}3*E03*E10\^{}2*E20+33*E02\^{}3\newline
*E10\^{}2*E20+14*E01\^{}2*E02\^{}2*E10\^{}2*E20-13*E01\^{}4*E02*E10\^{}2*E20+2*E01\^{}6\newline
*E10\^{}2*E20-27*E03*E10\^{}2*E11*E12+9*E01*E02*E10\^{}2*E11*E12-2*E01\^{}3\newline
*E10\^{}2*E11*E12+18*E01*E03*E10\^{}3*E12+6*E02\^{}2*E10\^{}3*E12-10*E01\^{}2*E02\newline
*E10\^{}3*E12+2*E01\^{}4*E10\^{}3*E12-9*E02*E11\^{}4+3*E01\^{}2*E11\^{}4+24*E01*E02\newline
*E10*E11\^{}3-8*E01\^{}3*E10*E11\^{}3+9*E01*E03*E10\^{}2*E11\^{}2-15*E02\^{}2*E10\^{}2\newline
*E11\^{}2-17*E01\^{}2*E02*E10\^{}2*E11\^{}2+7*E01\^{}4*E10\^{}2*E11\^{}2-9*E02*E03*E10\^{}3\newline
*E11-6*E01\^{}2*E03*E10\^{}3*E11+21*E01*E02\^{}2*E10\^{}3*E11-2*E01\^{}5*E10\^{}3*E11\newline
-27*E03\^{}2*E10\^{}4+24*E01*E02*E03*E10\^{}4-4*E01\^{}3*E03*E10\^{}4-6*E02\^{}3*E10\^{}4\newline
-5*E01\^{}2*E02\^{}2*E10\^{}4+2*E01\^{}4*E02*E10\^{}4;}
\medskip\noindent
{\tt\~q12:=\relax
-27*E02\^{}3*E30\^{}2+27*E01\^{}2*E02\^{}2*E30\^{}2-9*E01\^{}4*E02*E30\^{}2+E01\^{}6\newline
*E30\^{}2-81*E02*E03*E11*E20*E30+27*E01\^{}2*E03*E11*E20*E30+27*E01*E02\^{}2\newline
*E11*E20*E30-15*E01\^{}3*E02*E11*E20*E30+2*E01\^{}5*E11*E20*E30+54*E01*E02\newline
*E03*E10*E20*E30-18*E01\^{}3*E03*E10*E20*E30+18*E02\^{}3*E10*E20*E30-36\newline
*E01\^{}2*E02\^{}2*E10*E20*E30+16*E01\^{}4*E02*E10*E20*E30-2*E01\^{}6*E10*E20\newline
*E30+27*E03*E11\^{}3*E30-9*E01*E02*E11\^{}3*E30+2*E01\^{}3*E11\^{}3*E30-54*E01\newline
*E03*E10*E11\^{}2*E30+18*E01\^{}2*E02*E10*E11\^{}2*E30-4*E01\^{}4*E10*E11\^{}2*E30\newline
+27*E02*E03*E10\^{}2*E11*E30+27*E01\^{}2*E03*E10\^{}2*E11*E30-9*E01*E02\^{}2\newline
*E10\^{}2*E11*E30-7*E01\^{}3*E02*E10\^{}2*E11*E30+2*E01\^{}5*E10\^{}2*E11*E30-18\newline
*E01*E02*E03*E10\^{}3*E30-2*E01\^{}3*E03*E10\^{}3*E30-4*E02\^{}3*E10\^{}3*E30+10\newline
*E01\^{}2*E02\^{}2*E10\^{}3*E30-2*E01\^{}4*E02*E10\^{}3*E30-27*E03\^{}2*E20\^{}3+18*E01\newline
*E02*E03*E20\^{}3-4*E01\^{}3*E03*E20\^{}3-4*E02\^{}3*E20\^{}3+E01\^{}2*E02\^{}2*E20\^{}3+9\newline
*E02\^{}2*E11\^{}2*E20\^{}2-6*E01\^{}2*E02*E11\^{}2*E20\^{}2+E01\^{}4*E11\^{}2*E20\^{}2+27*E02\newline
*E03*E10*E11*E20\^{}2-9*E01\^{}2*E03*E10*E11*E20\^{}2-21*E01*E02\^{}2*E10*E11\newline
*E20\^{}2+13*E01\^{}3*E02*E10*E11*E20\^{}2-2*E01\^{}5*E10*E11*E20\^{}2+27*E03\^{}2\newline
*E10\^{}2*E20\^{}2-36*E01*E02*E03*E10\^{}2*E20\^{}2+10*E01\^{}3*E03*E10\^{}2*E20\^{}2\newline
+E02\^{}3*E10\^{}2*E20\^{}2+12*E01\^{}2*E02\^{}2*E10\^{}2*E20\^{}2-7*E01\^{}4*E02*E10\^{}2\newline
*E20\^{}2+E01\^{}6*E10\^{}2*E20\^{}2-6*E02*E11\^{}4*E20+2*E01\^{}2*E11\^{}4*E20-9*E03*E10\newline
*E11\^{}3*E20+19*E01*E02*E10*E11\^{}3*E20-6*E01\^{}3*E10*E11\^{}3*E20+18*E01*E03\newline
*E10\^{}2*E11\^{}2*E20-6*E02\^{}2*E10\^{}2*E11\^{}2*E20-18*E01\^{}2*E02*E10\^{}2*E11\^{}2\newline
*E20+6*E01\^{}4*E10\^{}2*E11\^{}2*E20-15*E02*E03*E10\^{}3*E11*E20-7*E01\^{}2*E03\newline
*E10\^{}3*E11*E20+13*E01*E02\^{}2*E10\^{}3*E11*E20+3*E01\^{}3*E02*E10\^{}3*E11*E20\newline
-2*E01\^{}5*E10\^{}3*E11*E20-9*E03\^{}2*E10\^{}4*E20+16*E01*E02*E03*E10\^{}4*E20-2\newline
*E01\^{}3*E03*E10\^{}4*E20-7*E01\^{}2*E02\^{}2*E10\^{}4*E20+2*E01\^{}4*E02*E10\^{}4*E20\newline
+E11\^{}6-4*E01*E10*E11\^{}5+2*E02*E10\^{}2*E11\^{}4+6*E01\^{}2*E10\^{}2*E11\^{}4+2*E03\newline
*E10\^{}3*E11\^{}3-6*E01*E02*E10\^{}3*E11\^{}3-4*E01\^{}3*E10\^{}3*E11\^{}3-4*E01*E03\newline
*E10\^{}4*E11\^{}2+E02\^{}2*E10\^{}4*E11\^{}2+6*E01\^{}2*E02*E10\^{}4*E11\^{}2+E01\^{}4*E10\^{}4\newline
*E11\^{}2+2*E02*E03*E10\^{}5*E11+2*E01\^{}2*E03*E10\^{}5*E11-2*E01*E02\^{}2*E10\^{}5\newline
*E11-2*E01\^{}3*E02*E10\^{}5*E11+E03\^{}2*E10\^{}6-2*E01*E02*E03*E10\^{}6+E01\^{}2\newline
*E02\^{}2*E10\^{}6;}
\medskip\noindent
{\tt\~q13:=\relax
-729*E02\^{}3*E12*E30\^{}2+729*E01\^{}2*E02\^{}2*E12*E30\^{}2-243*E01\^{}4*E02\newline
*E12*E30\^{}2+27*E01\^{}6*E12*E30\^{}2+2187*E02\^{}2*E03*E11*E30\^{}2-1458*E01\^{}2\newline
*E02*E03*E11*E30\^{}2+243*E01\^{}4*E03*E11*E30\^{}2-486*E01*E02\^{}3*E11*E30\^{}2\newline
+405*E01\^{}3*E02\^{}2*E11*E30\^{}2-108*E01\^{}5*E02*E11*E30\^{}2+9*E01\^{}7*E11*E30\^{}2\newline
-1458*E01*E02\^{}2*E03*E10*E30\^{}2+972*E01\^{}3*E02*E03*E10*E30\^{}2-162*E01\^{}5\newline
*E03*E10*E30\^{}2+243*E02\^{}4*E10*E30\^{}2+81*E01\^{}2*E02\^{}3*E10*E30\^{}2-189\newline
*E01\^{}4*E02\^{}2*E10*E30\^{}2+63*E01\^{}6*E02*E10*E30\^{}2-6*E01\^{}8*E10*E30\^{}2+4374\newline
*E02*E03\^{}2*E20\^{}2*E30-1458*E01\^{}2*E03\^{}2*E20\^{}2*E30-2916*E01*E02\^{}2*E03\newline
*E20\^{}2*E30+1620*E01\^{}3*E02*E03*E20\^{}2*E30-216*E01\^{}5*E03*E20\^{}2*E30+324\newline
*E02\^{}4*E20\^{}2*E30+54*E01\^{}2*E02\^{}3*E20\^{}2*E30-162*E01\^{}4*E02\^{}2*E20\^{}2*E30\newline
+48*E01\^{}6*E02*E20\^{}2*E30-4*E01\^{}8*E20\^{}2*E30+486*E02\^{}3*E10*E12*E20*E30\newline
-486*E01\^{}2*E02\^{}2*E10*E12*E20*E30+162*E01\^{}4*E02*E10*E12*E20*E30-18\newline
*E01\^{}6*E10*E12*E20*E30-2187*E03\^{}2*E11\^{}2*E20*E30+1458*E01*E02*E03\newline
*E11\^{}2*E20*E30-324*E01\^{}3*E03*E11\^{}2*E20*E30-405*E02\^{}3*E11\^{}2*E20*E30\newline
+162*E01\^{}2*E02\^{}2*E11\^{}2*E20*E30-27*E01\^{}4*E02*E11\^{}2*E20*E30+3*E01\^{}6\newline
*E11\^{}2*E20*E30+2916*E01*E03\^{}2*E10*E11*E20*E30-1458*E02\^{}2*E03*E10*E11\newline
*E20*E30-972*E01\^{}2*E02*E03*E10*E11*E20*E30+270*E01\^{}4*E03*E10*E11*E20\newline
*E30+864*E01*E02\^{}3*E10*E11*E20*E30-486*E01\^{}3*E02\^{}2*E10*E11*E20*E30\newline
+108*E01\^{}5*E02*E10*E11*E20*E30-10*E01\^{}7*E10*E11*E20*E30-2916*E02\newline
*E03\^{}2*E10\^{}2*E20*E30+2916*E01*E02\^{}2*E03*E10\^{}2*E20*E30-1080*E01\^{}3*E02\newline
*E03*E10\^{}2*E20*E30+108*E01\^{}5*E03*E10\^{}2*E20*E30-378*E02\^{}4*E10\^{}2*E20\newline
*E30-270*E01\^{}2*E02\^{}3*E10\^{}2*E20*E30+306*E01\^{}4*E02\^{}2*E10\^{}2*E20*E30-86\newline
*E01\^{}6*E02*E10\^{}2*E20*E30+8*E01\^{}8*E10\^{}2*E20*E30-108*E02\^{}3*E10\^{}3*E12\newline
*E30+108*E01\^{}2*E02\^{}2*E10\^{}3*E12*E30-36*E01\^{}4*E02*E10\^{}3*E12*E30+4\newline
*E01\^{}6*E10\^{}3*E12*E30+81*E02\^{}2*E11\^{}4*E30-54*E01\^{}2*E02*E11\^{}4*E30+9\newline
*E01\^{}4*E11\^{}4*E30-216*E01*E02\^{}2*E10*E11\^{}3*E30+144*E01\^{}3*E02*E10*E11\^{}3\newline
*E30-24*E01\^{}5*E10*E11\^{}3*E30+729*E03\^{}2*E10\^{}2*E11\^{}2*E30-486*E01*E02\newline
*E03*E10\^{}2*E11\^{}2*E30+108*E01\^{}3*E03*E10\^{}2*E11\^{}2*E30+135*E02\^{}3*E10\^{}2\newline
*E11\^{}2*E30+162*E01\^{}2*E02\^{}2*E10\^{}2*E11\^{}2*E30-135*E01\^{}4*E02*E10\^{}2*E11\^{}2\newline
*E30+23*E01\^{}6*E10\^{}2*E11\^{}2*E30-972*E01*E03\^{}2*E10\^{}3*E11*E30+324*E02\^{}2\newline
*E03*E10\^{}3*E11*E30+432*E01\^{}2*E02*E03*E10\^{}3*E11*E30-108*E01\^{}4*E03\newline
*E10\^{}3*E11*E30-252*E01*E02\^{}3*E10\^{}3*E11*E30+36*E01\^{}3*E02\^{}2*E10\^{}3*E11\newline
*E30+36*E01\^{}5*E02*E10\^{}3*E11*E30-8*E01\^{}7*E10\^{}3*E11*E30+486*E02*E03\^{}2\newline
*E10\^{}4*E30+162*E01\^{}2*E03\^{}2*E10\^{}4*E30-540*E01*E02\^{}2*E03*E10\^{}4*E30+108\newline
*E01\^{}3*E02*E03*E10\^{}4*E30+72*E02\^{}4*E10\^{}4*E30+78*E01\^{}2*E02\^{}3*E10\^{}4*E30\newline
-54*E01\^{}4*E02\^{}2*E10\^{}4*E30+8*E01\^{}6*E02*E10\^{}4*E30+729*E03\^{}2*E12*E20\^{}3\newline
\vadjust{\vskip 0pt plus 1pt minus 1 pt}\relax
-486*E01*E02*E03*E12*E20\^{}3+108*E01\^{}3*E03*E12*E20\^{}3+81*E01\^{}2*E02\^{}2\newline
*E12*E20\^{}3-36*E01\^{}4*E02*E12*E20\^{}3+4*E01\^{}6*E12*E20\^{}3-243*E01*E03\^{}2\newline
*E11*E20\^{}3-324*E02\^{}2*E03*E11*E20\^{}3+378*E01\^{}2*E02*E03*E11*E20\^{}3-72\newline
*E01\^{}4*E03*E11*E20\^{}3+108*E01*E02\^{}3*E11*E20\^{}3-123*E01\^{}3*E02\^{}2*E11\newline
*E20\^{}3+40*E01\^{}5*E02*E11*E20\^{}3-4*E01\^{}7*E11*E20\^{}3-1701*E02*E03\^{}2*E10\newline
*E20\^{}3+648*E01\^{}2*E03\^{}2*E10*E20\^{}3+1350*E01*E02\^{}2*E03*E10*E20\^{}3-828\newline
*E01\^{}3*E02*E03*E10*E20\^{}3+120*E01\^{}5*E03*E10*E20\^{}3-108*E02\^{}4*E10*E20\^{}3\newline
-117*E01\^{}2*E02\^{}3*E10*E20\^{}3+148*E01\^{}4*E02\^{}2*E10*E20\^{}3-44*E01\^{}6*E02\newline
*E10*E20\^{}3+4*E01\^{}8*E10*E20\^{}3-729*E03\^{}2*E10\^{}2*E12*E20\^{}2+486*E01*E02\newline
*E03*E10\^{}2*E12*E20\^{}2-108*E01\^{}3*E03*E10\^{}2*E12*E20\^{}2-81*E02\^{}3*E10\^{}2\newline
*E12*E20\^{}2+9*E01\^{}4*E02*E10\^{}2*E12*E20\^{}2-E01\^{}6*E10\^{}2*E12*E20\^{}2+405*E02\newline
*E03*E11\^{}3*E20\^{}2-135*E01\^{}2*E03*E11\^{}3*E20\^{}2-135*E01*E02\^{}2*E11\^{}3*E20\^{}2\newline
+75*E01\^{}3*E02*E11\^{}3*E20\^{}2-10*E01\^{}5*E11\^{}3*E20\^{}2+729*E03\^{}2*E10*E11\^{}2\newline
*E20\^{}2-1296*E01*E02*E03*E10*E11\^{}2*E20\^{}2+378*E01\^{}3*E03*E10*E11\^{}2\newline
*E20\^{}2+135*E02\^{}3*E10*E11\^{}2*E20\^{}2+216*E01\^{}2*E02\^{}2*E10*E11\^{}2*E20\^{}2-141\newline
*E01\^{}4*E02*E10*E11\^{}2*E20\^{}2+19*E01\^{}6*E10*E11\^{}2*E20\^{}2-729*E01*E03\^{}2\newline
*E10\^{}2*E11*E20\^{}2+567*E02\^{}2*E03*E10\^{}2*E11*E20\^{}2+648*E01\^{}2*E02*E03\newline
*E10\^{}2*E11*E20\^{}2-225*E01\^{}4*E03*E10\^{}2*E11*E20\^{}2-342*E01*E02\^{}3*E10\^{}2\newline
*E11*E20\^{}2+60*E01\^{}3*E02\^{}2*E10\^{}2*E11*E20\^{}2+36*E01\^{}5*E02*E10\^{}2*E11\newline
*E20\^{}2-7*E01\^{}7*E10\^{}2*E11*E20\^{}2+1539*E02*E03\^{}2*E10\^{}3*E20\^{}2-270*E01\^{}2\newline
*E03\^{}2*E10\^{}3*E20\^{}2-1404*E01*E02\^{}2*E03*E10\^{}3*E20\^{}2+540*E01\^{}3*E02*E03\newline
*E10\^{}3*E20\^{}2-42*E01\^{}5*E03*E10\^{}3*E20\^{}2+123*E02\^{}4*E10\^{}3*E20\^{}2+184\newline
*E01\^{}2*E02\^{}3*E10\^{}3*E20\^{}2-147*E01\^{}4*E02\^{}2*E10\^{}3*E20\^{}2+31*E01\^{}6*E02\newline
*E10\^{}3*E20\^{}2-2*E01\^{}8*E10\^{}3*E20\^{}2+243*E03\^{}2*E10\^{}4*E12*E20-162*E01*E02\newline
*E03*E10\^{}4*E12*E20+36*E01\^{}3*E03*E10\^{}4*E12*E20+36*E02\^{}3*E10\^{}4*E12*E20\newline
-9*E01\^{}2*E02\^{}2*E10\^{}4*E12*E20-81*E03*E11\^{}5*E20+27*E01*E02*E11\^{}5*E20-6\newline
*E01\^{}3*E11\^{}5*E20+270*E01*E03*E10*E11\^{}4*E20-27*E02\^{}2*E10*E11\^{}4*E20-72\newline
*E01\^{}2*E02*E10*E11\^{}4*E20+17*E01\^{}4*E10*E11\^{}4*E20-270*E02*E03*E10\^{}2\newline
*E11\^{}3*E20-270*E01\^{}2*E03*E10\^{}2*E11\^{}3*E20+162*E01*E02\^{}2*E10\^{}2*E11\^{}3\newline
*E20+22*E01\^{}3*E02*E10\^{}2*E11\^{}3*E20-12*E01\^{}5*E10\^{}2*E11\^{}3*E20-405*E03\^{}2\newline
*E10\^{}3*E11\^{}2*E20+810*E01*E02*E03*E10\^{}3*E11\^{}2*E20-75*E02\^{}3*E10\^{}3\newline
*E11\^{}2*E20-222*E01\^{}2*E02\^{}2*E10\^{}3*E11\^{}2*E20+63*E01\^{}4*E02*E10\^{}3*E11\^{}2\newline
*E20-3*E01\^{}6*E10\^{}3*E11\^{}2*E20+459*E01*E03\^{}2*E10\^{}4*E11*E20-216*E02\^{}2\newline
*E03*E10\^{}4*E11*E20-522*E01\^{}2*E02*E03*E10\^{}4*E11*E20+84*E01\^{}4*E03\newline
*E10\^{}4*E11*E20+160*E01*E02\^{}3*E10\^{}4*E11*E20+51*E01\^{}3*E02\^{}2*E10\^{}4*E11\newline
*E20-36*E01\^{}5*E02*E10\^{}4*E11*E20+4*E01\^{}7*E10\^{}4*E11*E20-459*E02*E03\^{}2\newline
*E10\^{}5*E20+450*E01*E02\^{}2*E03*E10\^{}5*E20-84*E01\^{}3*E02*E03*E10\^{}5*E20-40\newline
*E02\^{}4*E10\^{}5*E20-75*E01\^{}2*E02\^{}3*E10\^{}5*E20+36*E01\^{}4*E02\^{}2*E10\^{}5*E20-4\newline
*E01\^{}6*E02*E10\^{}5*E20-27*E03\^{}2*E10\^{}6*E12+18*E01*E02*E03*E10\^{}6*E12-4\newline
*E01\^{}3*E03*E10\^{}6*E12-4*E02\^{}3*E10\^{}6*E12+E01\^{}2*E02\^{}2*E10\^{}6*E12+27*E03\newline
*E10\^{}2*E11\^{}5-9*E01*E02*E10\^{}2*E11\^{}5+2*E01\^{}3*E10\^{}2*E11\^{}5-90*E01*E03\newline
*E10\^{}3*E11\^{}4+6*E02\^{}2*E10\^{}3*E11\^{}4+26*E01\^{}2*E02*E10\^{}3*E11\^{}4-6*E01\^{}4\newline
*E10\^{}3*E11\^{}4+45*E02*E03*E10\^{}4*E11\^{}3+105*E01\^{}2*E03*E10\^{}4*E11\^{}3-31*E01\newline
*E02\^{}2*E10\^{}4*E11\^{}3-21*E01\^{}3*E02*E10\^{}4*E11\^{}3+6*E01\^{}5*E10\^{}4*E11\^{}3+54\newline
*E03\^{}2*E10\^{}5*E11\^{}2-126*E01*E02*E03*E10\^{}5*E11\^{}2-42*E01\^{}3*E03*E10\^{}5\newline
*E11\^{}2+10*E02\^{}3*E10\^{}5*E11\^{}2+42*E01\^{}2*E02\^{}2*E10\^{}5*E11\^{}2-2*E01\^{}6*E10\^{}5\newline
*E11\^{}2-63*E01*E03\^{}2*E10\^{}6*E11+24*E02\^{}2*E03*E10\^{}6*E11+86*E01\^{}2*E02\newline
*E03*E10\^{}6*E11-20*E01*E02\^{}3*E10\^{}6*E11-15*E01\^{}3*E02\^{}2*E10\^{}6*E11+4\newline
*E01\^{}5*E02*E10\^{}6*E11+45*E02*E03\^{}2*E10\^{}7+6*E01\^{}2*E03\^{}2*E10\^{}7-46*E01\newline
*E02\^{}2*E03*E10\^{}7+4*E02\^{}4*E10\^{}7+9*E01\^{}2*E02\^{}3*E10\^{}7-2*E01\^{}4*E02\^{}2\newline
*E10\^{}7;}
\medskip\noindent
{\tt\~q14:=\relax
81*E03\^{}2*E30-27*E01*E02*E03*E30+3*E01\^{}3*E03*E30+3*E02\^{}3*E30+9\newline
*E01*E03*E12*E20-3*E02\^{}2*E12*E20-9*E02*E03*E11*E20+3*E01\^{}2*E03*E11\newline
*E20-54*E03\^{}2*E10*E20+30*E01*E02*E03*E10*E20-7*E01\^{}3*E03*E10*E20-4\newline
*E02\^{}3*E10*E20+E01\^{}2*E02\^{}2*E10*E20-3*E12\^{}3+E01\^{}2*E10*E12\^{}2+9*E03*E10\newline
*E11*E12-E01*E02*E10*E11*E12-9*E01*E03*E10\^{}2*E12+2*E02\^{}2*E10\^{}2*E12+3\newline
*E03*E11\^{}3-9*E01*E03*E10*E11\^{}2+E02\^{}2*E10*E11\^{}2+7*E01\^{}2*E03*E10\^{}2*E11\newline
-E01*E02\^{}2*E10\^{}2*E11+15*E03\^{}2*E10\^{}3-7*E01*E02*E03*E10\^{}3+E02\^{}3*E10\^{}3.}
\par

\Refs
\ref\myrefno{1}\paper
\myhref{http://en.wikipedia.org/wiki/Euler\podcherkivanie 
brick}{Euler brick}\jour Wikipedia\publ 
Wikimedia Foundation Inc.\publaddr San Francisco, USA 
\endref
\ref\myrefno{2}\by Halcke~P.\book Deliciae mathematicae oder mathematisches 
Sinnen-Confect\publ N.~Sauer\publaddr Hamburg, Germany\yr 1719
\endref
\ref\myrefno{3}\by Saunderson~N.\book Elements of algebra, {\rm Vol. 2}\publ
Cambridge Univ\. Press\publaddr Cambridge\yr 1740 
\endref
\ref\myrefno{4}\by Euler~L.\book Vollst\"andige Anleitung zur Algebra, \rm
3 Theile\publ Kaiserliche Akademie der Wissenschaf\-ten\publaddr St\.~Petersburg
\yr 1770-1771
\endref
\ref\myrefno{5}\by Pocklington~H.~C.\paper Some Diophantine impossibilities
\jour Proc. Cambridge Phil\. Soc\. \vol 17\yr 1912\pages 108--121
\endref
\ref\myrefno{6}\by Dickson~L.~E\book History of the theory of numbers, 
{\rm Vol\. 2}: Diophantine analysis\publ Dover\publaddr New York\yr 2005
\endref
\ref\myrefno{7}\by Kraitchik~M.\paper On certain rational cuboids
\jour Scripta Math\.\vol 11\yr 1945\pages 317--326
\endref
\ref\myrefno{8}\by Kraitchik~M.\book Th\'eorie des Nombres,
{\rm Tome 3}, Analyse Diophantine et application aux cuboides 
rationelles \publ Gauthier-Villars\publaddr Paris\yr 1947
\endref
\ref\myrefno{9}\by Kraitchik~M.\paper Sur les cuboides rationelles
\jour Proc\. Int\. Congr\. Math\.\vol 2\yr 1954\publaddr Amsterdam
\pages 33--34
\endref
\ref\myrefno{10}\by Bromhead~T.~B.\paper On square sums of squares
\jour Math\. Gazette\vol 44\issue 349\yr 1960\pages 219--220
\endref
\ref\myrefno{11}\by Lal~M., Blundon~W.~J.\paper Solutions of the 
Diophantine equations $x^2+y^2=l^2$, $y^2+z^2=m^2$, $z^2+x^2
=n^2$\jour Math\. Comp\.\vol 20\yr 1966\pages 144--147
\endref
\ref\myrefno{12}\by Spohn~W.~G.\paper On the integral cuboid\jour Amer\. 
Math\. Monthly\vol 79\issue 1\pages 57-59\yr 1972 
\endref
\ref\myrefno{13}\by Spohn~W.~G.\paper On the derived cuboid\jour Canad\. 
Math\. Bull\.\vol 17\issue 4\pages 575-577\yr 1974
\endref
\ref\myrefno{14}\by Chein~E.~Z.\paper On the derived cuboid of an 
Eulerian triple\jour Canad\. Math\. Bull\.\vol 20\issue 4\yr 1977
\pages 509--510
\endref
\ref\myrefno{15}\by Leech~J.\paper The rational cuboid revisited
\jour Amer\. Math\. Monthly\vol 84\issue 7\pages 518--533\yr 1977
\moreref see also Erratum\jour Amer\. Math\. Monthly\vol 85\page 472
\yr 1978
\endref
\ref\myrefno{16}\by Leech~J.\paper Five tables relating to rational cuboids
\jour Math\. Comp\.\vol 32\yr 1978\pages 657--659
\endref
\ref\myrefno{17}\by Spohn~W.~G.\paper Table of integral cuboids and their 
generators\jour Math\. Comp\.\vol 33\yr 1979\pages 428--429
\endref
\ref\myrefno{18}\by Lagrange~J.\paper Sur le d\'eriv\'e du cuboide 
Eul\'erien\jour Canad\. Math\. Bull\.\vol 22\issue 2\yr 1979\pages 239--241
\endref
\ref\myrefno{19}\by Leech~J.\paper A remark on rational cuboids\jour Canad\. 
Math\. Bull\.\vol 24\issue 3\yr 1981\pages 377--378
\endref
\ref\myrefno{20}\by Korec~I.\paper Nonexistence of small perfect 
rational cuboid\jour Acta Math\. Univ\. Comen\.\vol 42/43\yr 1983
\pages 73--86
\endref
\ref\myrefno{21}\by Korec~I.\paper Nonexistence of small perfect 
rational cuboid II\jour Acta Math\. Univ\. Comen\.\vol 44/45\yr 1984
\pages 39--48
\endref
\ref\myrefno{22}\by Wells~D.~G.\book The Penguin dictionary of curious and 
interesting numbers\publ Penguin publishers\publaddr London\yr 1986
\endref
\ref\myrefno{23}\by Bremner~A., Guy~R.~K.\paper A dozen difficult Diophantine 
dilemmas\jour Amer\. Math\. Monthly\vol 95\issue 1\yr 1988\pages 31--36
\endref
\ref\myrefno{24}\by Bremner~A.\paper The rational cuboid and a quartic surface
\jour Rocky Mountain J\. Math\. \vol 18\issue 1\yr 1988\pages 105--121
\endref
\ref\myrefno{25}\by Colman~W.~J.~A.\paper On certain semiperfect cuboids\jour
Fibonacci Quart.\vol 26\issue 1\yr 1988\pages 54--57\moreref see also\nofrills 
\paper Some observations on the classical cuboid and its parametric solutions
\jour Fibonacci Quart\.\vol 26\issue 4\yr 1988\pages 338--343
\endref
\ref\myrefno{26}\by Korec~I.\paper Lower bounds for perfect rational cuboids 
\jour Math\. Slovaca\vol 42\issue 5\yr 1992\pages 565--582
\endref
\ref\myrefno{27}\by Guy~R.~K.\paper Is there a perfect cuboid? Four squares 
whose sums in pairs are square. Four squares whose differences are square 
\inbook Unsolved Problems in Number Theory, 2nd ed.\pages 173--181\yr 1994
\publ Springer-Verlag\publaddr New York 
\endref
\ref\myrefno{28}\by Rathbun~R.~L., Granlund~T.\paper The integer cuboid table 
with body, edge, and face type of solutions\jour Math\. Comp\.\vol 62\yr 1994
\pages 441--442
\endref
\ref\myrefno{29}\by Van Luijk~R.\book On perfect cuboids, \rm Doctoraalscriptie
\publ Mathematisch Instituut, Universiteit Utrecht\publaddr Utrecht\yr 2000
\endref
\ref\myrefno{30}\by Rathbun~R.~L., Granlund~T.\paper The classical rational 
cuboid table of Maurice Kraitchik\jour Math\. Comp\.\vol 62\yr 1994
\pages 442--443
\endref
\ref\myrefno{31}\by Peterson~B.~E., Jordan~J.~H.\paper Integer hexahedra equivalent 
to perfect boxes\jour Amer\. Math\. Monthly\vol 102\issue 1\yr 1995\pages 41--45
\endref
\ref\myrefno{32}\by Rathbun~R.~L.\paper The rational cuboid table of Maurice 
Kraitchik\jour e-print \myhref{http://arxiv.org/abs/math/0111229}{math.HO/0111229} 
in Electronic Archive \myEarXivlink
\endref
\ref\myrefno{33}\by Hartshorne~R., Van Luijk~R.\paper Non-Euclidean Pythagorean 
triples, a problem of Euler, and rational points on K3 surfaces\publ e-print 
\myhref{http://arxiv.org/abs/math/0606700}{math.NT/0606700} 
in Electronic Archive \myEarXivlink
\endref
\ref\myrefno{34}\by Waldschmidt~M.\paper Open diophantine problems\publ e-print 
\myhref{http://arxiv.org/abs/math/0312440}{math.NT/0312440} 
in Electronic Archive \myEarXivlink
\endref
\ref\myrefno{35}\by Ionascu~E.~J., Luca~F., Stanica~P.\paper Heron triangles 
with two fixed sides\publ e-print \myhref{http://arxiv.org/abs/math/0608185}
{math.NT/0608} \myhref{http://arxiv.org/abs/math/0608185}{185} in Electronic 
Archive \myEarXivlink
\endref
\ref\myrefno{36}\by Ortan~A., Quenneville-Belair~V.\paper Euler's brick
\jour Delta Epsilon, McGill Undergraduate Mathematics Journal\yr 2006\vol 1
\pages 30-33
\endref
\ref\myrefno{37}\by Knill~O.\paper Hunting for Perfect Euler Bricks\jour Harvard
College Math\. Review\yr 2008\vol 2\issue 2\page 102\moreref
see also \myhref{http://www.math.harvard.edu/\volna knill/various/eulercuboid/index.html}
{http:/\negskp/www.math.harvard.edu/\textvolna knill/various/eulercuboid/index.html}
\endref
\ref\myrefno{38}\by Sloan~N.~J.~A\paper Sequences 
\myhref{http://oeis.org/A031173}{A031173}, 
\myhref{http://oeis.org/A031174}{A031174}, and \myhref{http://oeis.org/A031175}
{A031175}\jour On-line encyclopedia of integer sequences\publ OEIS Foundation 
Inc.\publaddr Portland, USA
\endref
\ref\myrefno{39}\by Stoll~M., Testa~D.\paper The surface parametrizing cuboids
\jour e-print \myhref{http://arxiv.org/abs/1009.0388}{arXiv:1009.0388} 
in Electronic Archive \myEarXivlink
\endref
\ref\myrefno{40}\by Sharipov~R.~A.\paper A note on a perfect Euler cuboid.
\jour e-print \myhref{http://arxiv.org/abs/1104.1716}{arXiv:1104.1716} 
in Electronic Archive \myEarXivlink
\endref
\ref\myrefno{41}\by Sharipov~R.~A.\paper Perfect cuboids and irreducible 
polynomials\jour Ufa Mathematical Journal\vol 4, \issue 1\yr 2012\pages 153--160
\moreref see also e-print \myhref{http://arxiv.org/abs/1108.5348}{arXiv:1108.5348} 
in Electronic Archive \myEarXivlink
\endref
\ref\myrefno{42}\by Sharipov~R.~A.\paper A note on the first cuboid conjecture
\jour e-print \myhref{http://arxiv.org/abs/1109.2534}{arXiv:1109.2534} 
in Electronic Archive \myEarXivlink
\endref
\ref\myrefno{43}\by Sharipov~R.~A.\paper A note on the second cuboid conjecture.
Part~\uppercase\expandafter{\romannumeral 1} 
\jour e-print \myhref{http://arxiv.org/abs/1201.1229}{arXiv:1201.1229} 
in Electronic Archive \myEarXivlink
\endref
\ref\myrefno{44}\by Sharipov~R.~A.\paper A note on the third cuboid conjecture.
Part~\uppercase\expandafter{\romannumeral 1} 
\jour e-print \myhref{http://arxiv.org/abs/1203.2567}{arXiv:1203.2567} 
in Electronic Archive \myEarXivlink
\endref
\ref\myrefno{45}\by Sharipov~R.~A.\paper Perfect cuboids and multisymmetric 
polynomials\jour e-print \myhref{http://arxiv.org/abs/1203.2567}
{arXiv:1205.3135} in Electronic Archive \myEarXivlink
\endref
\ref\myrefno{46}\by Shl\"afli~L.\paper \"Uber die Resultante eines systems mehrerer 
algebraishen Gleihungen\jour Denkschr\. Kaiserliche Acad\. Wiss\. Math\.-Natur\.
Kl\.\vol 4\yr 1852\moreref reprinted in {\eightcyr\char '074}Gesammelte mathematische
Abhandlungen{\eightcyr\char '076}, Band \uppercase\expandafter{\romannumeral 2}
\pages 9--112\publ Birkh\"auser Verlag\yr 1953
\endref
\ref\myrefno{47}\by Cayley~A.\paper On the symmetric functions of the roots of 
certain systems of two equations\jour Phil\. Trans\. Royal Soc\. London\vol 147
\yr 1857\pages 717--726
\endref
\ref\myrefno{48}\by Junker~F.\paper \"Uber symmetrische Functionen von mehreren 
Ver\"anderlishen\jour Mathematische Annalen\vol 43\pages 225--270 \yr 1893
\endref
\ref\myrefno{49}\by McMahon~P.~A.\paper Memoir on symmetric functions of the
roots of systems of equations\jour Phil\. Trans\. Royal Soc\. London\vol 181
\yr 1890\pages 481--536
\endref
\ref\myrefno{50}\by McMahon~P.~A. \book Combinatory Analysis. 
\rm Vol\.~\uppercase\expandafter{\romannumeral 1} and 
Vol\.~\uppercase\expandafter{\romannumeral 2}\publ Cambridge Univ\. Press
\yr 1915--1916\moreref see also Third ed\.\publ Chelsea Publishing Company
\publaddr New York\yr 1984
\endref
\ref\myrefno{51}\by Noether~E.\paper Der Endlichkeitssats der Invarianten
endlicher Gruppen\jour Mathematische Annalen\vol 77\pages 89--92 \yr 1915
\endref
\ref\myrefno{52}\by Weyl~H.\book The classical groups\publ Princeton Univ\.
Press\publaddr Princeton\yr1939
\endref
\ref\myrefno{53}\by Macdonald~I.~G.\book Symmetric functions and Hall polynomials,
\rm Oxford Mathematical Monographs\publ Clarendon Press\publaddr Oxford\yr 1979 
\endref
\ref\myrefno{54}\by Pedersen~P.\paper Calculating multidimensional symmetric
functions using Jacobi's formula\inbook Proceedings AAECC 9, volume 539 of
Springer Lecture Notes in Computer Science\pages 304--317\yr 1991\publ Springer
\endref
\ref\myrefno{55}\by Milne~P.\paper On the solutions of a set of polynomial equations
\inbook Symbolic and numerical computation for artificial intelligence. Computational 
Mathematics and Applications\eds Donald~B.~R., Kapur~D., Mundy~J.~L.\yr 1992\publ
Academic Press Ltd.\publaddr London\pages 89--101
\endref
\ref\myrefno{56}\by Dalbec~J.\book Geometry and combinatorics of Chow forms
\publ PhD thesis, Cornell University\yr 1995
\endref
\ref\myrefno{57}\by Richman~D.~R.\paper Explicit generators of the invariants of 
finite groups\jour Advances in Math\.\vol 124\issue 1\yr 1996\pages 49--76
\endref
\ref\myrefno{58}\by Stepanov~S.~A.\paper On vector invariants of the symmetric group
\jour Diskretnaya Matematika\vol 8\issue 2\yr 1996\pages 48--62
\endref
\ref\myrefno{59}\by Gonzalez-Vega~L., Trujillo~G.\paper Multivariate Sturm-Habicht 
sequences: real root counting on n-rectangles and triangles\jour Revista Matem\'atica 
Complutense\vol 10\pages 119--130\yr 1997
\endref
\ref\myrefno{60}\by Stepanov~S.~A.\paper On vector invariants of symmetric groups
\jour Diskretnaya Matematika\vol 11\issue 3\yr 1999\pages 4--14
\endref
\ref\myrefno{61}\by Dalbec~J.\paper Multisymmetric functions\jour Beitr\"age zur
Algebra und Geom\.\vol 40\issue 1\yr 1999\pages 27--51
\endref
\ref\myrefno{62}\by Rosas~M.~H.\paper MacMahon symmetric functions, the partition 
lattice, and Young subgroups\jour Journ\. Combin. Theory\vol 96\,A\issue 2\yr 2001
\pages 326--340
\endref
\ref\myrefno{63}\by Vaccarino~F.\paper The ring of  multisymmetric functions
\jour e-print \myhref{http://arxiv.org/abs/math/0205233}{math.RA/0205233} 
in Electronic Archive \myEarXivlink
\endref
\ref\myrefno{64}\by Briand~E.\paper When is the algebra of multisymmetric 
polynomials generated by the elementary multisymmetric polynomials?
\jour Beitr\"age zur Algebra und Geom\.\vol 45 \issue 2\pages 353--368
\yr 2004
\endref
\ref\myrefno{65}\by Rota~G.-C., Stein~J.~A.\paper A problem of Cayley from 1857
and how he could have solved it\jour Linear Algebra and its Applications (special 
issue on determinants and the legacy of Sir Thomas Muir)\vol 411\pages 167--253
\yr 2005
\endref
\ref\myrefno{66}\by Briand~E., Rosas~M.~H.\paper Milne's volume function and vector 
symmetric polynomials\jour Journ. Symbolic Comput. \vol 44\issue 5\yr 2009
\pages 583--590
\endref
\ref\myrefno{67}\by Sharipov~R.~A.\paper On the equivalence of cuboid equations and 
their factor equations\jour e-print \myhref{http://arxiv.org/abs/1207.2102}
{arXiv:1207.2102} in Electronic Archive \myEarXivlink
\endref
\ref\myrefno{68}\by Sharipov~R.~A.\paper On an ideal of multisymmetric polynomials 
associated with perfect cuboids\jour e-print \myhref{http://arxiv.org/abs/1206.6769}
{arXiv:1206.6769} in Electronic Archive \myEarXivlink
\endref
\ref\myrefno{69}\paper \myhref{http://en.wikipedia.org/wiki/Symmetric\podcherkivanie
polynomial}{Symmetric polynomial}\jour Wikipedia\publ Wikimedia Foundation Inc.
\publaddr San Francisco, USA 
\endref
\ref\myrefno{70}\by Cox~D.~A., Little~J.~B., O'Shea~D.\book Ideals, Varieties, 
and Algorithms\publ Springer Verlag\publaddr New York\yr 1992
\endref
\ref\myrefno{71}\paper \myhref{http://en.wikipedia.org/wiki/Heronian\podcherkivanie
triangle}{Heronian triangle}\jour Wikipedia\publ Wikimedia Foundation Inc.
\publaddr San Francisco, USA 
\endref


\endRefs
\enddocument
\end